\documentclass[reqno,10pt,a4paper]{amsart} % no b5paper option

\usepackage{kotex}

\usepackage{amsthm}
\usepackage{amsmath}
\usepackage{amsfonts}
\usepackage{amssymb}

\usepackage{mathrsfs}% Support for using RSFS fonts in maths.  Provides \mtahscr overwriting \mathcal. It is necessary for the multiplier ideal command \MI.
\usepackage{stmaryrd} %St. Mary Road Symbol. for right(left)lifting property.
\usepackage{mathabx} % for \Asterisk
\usepackage[T1]{fontenc}% for bold small cap in section title. For example, \textsc{i}
\usepackage{dsfont}% for \mathds

\usepackage{graphicx}
\usepackage[all,2cell]{xy}
\usepackage{tikz}
\usetikzlibrary{matrix,arrows,decorations.pathmorphing,calc}
\makeatletter
\ifcsname tikz@library@external@loaded \endcsname
\else
\usepackage{tikz-cd}
\fi
\makeatother

\usepackage[colorlinks=true]{hyperref}
\usepackage{url}

\usepackage{setspace}% for line spacing

\usepackage{ifthen}%\ifthenelse{}{}{}

\usepackage{enumerate}

\usepackage{comment}

\usepackage{bm}

\usepackage{enumitem}

\setlistdepth{9}
\renewlist{itemize}{itemize}{9}
\setlist[itemize,1]{label=\textbullet}
\setlist[itemize,2]{label=--}
\setlist[itemize,3]{label=*}
\setlist[itemize,4]{label=\textbullet}
\setlist[itemize,5]{label=--}
\setlist[itemize,6]{label=*}
\setlist[itemize,7]{label=\textbullet}
\setlist[itemize,8]{label=--}
\setlist[itemize,9]{label=*}

\makeatletter
\@ifclassloaded{book}{
\usepackage[nottoc]{tocbibind}
}{}
\makeatother

\usepackage[normalem]{ulem}
\usepackage{cancel}

\makeatletter
\@ifclassloaded{book}{
  \usepackage{fancyhdr}

  \pagestyle{fancy}

  \fancyhf{}
  \fancyhead[LE,RO]{\bfseries\thepage}
  \fancyhead[LO]{\bfseries\rightmark}
  \fancyhead[RE]{\bfseries\leftmark}

  \addtolength{\headheight}{0.5pt} % space for the rule
  \fancypagestyle{plain}{%
    \fancyhead{} % get rid of headers on plain pages
  }
}{}
\makeatother

\makeatletter
\@ifclassloaded{book}{
  \usepackage{minitoc}
}{}
\makeatother

\makeatletter
\@ifclassloaded{book}{
  \usepackage{makeidx}

  \makeindex
}{}

\makeatletter
\@ifclassloaded{book}
{
  \setcounter{tocdepth}{3}%Max: 5
  \setcounter{secnumdepth}{5}
  \setcounter{minitocdepth}{5}
}
{
}

\makeatletter
\@ifclassloaded{book}{
  \newcounter{SaveAppChapter}
  \newcounter{SaveTextChapter}
  \edef\CurrentChapter{\thechapter}

  \newboolean{appchapter}
  \setboolean{appchapter}{false}

  \newcommand{\switchchapter}
  {
    \ifthenelse{
      \boolean{appchapter}
    }
    {
      \setcounter{SaveAppChapter}{\value{chapter}}
      \setboolean{appchapter}{false}
      \setcounter{chapter}{\value{SaveTextChapter}}
      \renewcommand{\thechapter}{\arabic{chapter}}
      \renewcommand{\theHchapter}{\arabic{chapter}}      
    }
    {
      \setcounter{SaveTextChapter}{\value{chapter}}
      \setboolean{appchapter}{true}
      \setcounter{chapter}{\value{SaveAppChapter}}
      \renewcommand{\thechapter}{\Alph{chapter}}
      \renewcommand{\theHchapter}{\Alph{chapter}}      
    }
  }
}{}

\newcounter{SaveAppSection}
\newcounter{SaveTextSection}
\edef\CurrentSection{\thesection}

\newboolean{appsection}
\setboolean{appsection}{false}

\newcommand{\switchsection}
{
  \ifthenelse{
    \boolean{appsection}
  }
  {
    \setcounter{SaveAppSection}{\value{section}}
    \setboolean{appsection}{false}
    \ifthenelse{\equal{\thechapter}{\CurrentChapter}}{}
    {
      \setcounter{SaveTextSection}{0}
      \edef\CurrentChapter{\thechapter}
    }
    \setcounter{section}{\value{SaveTextSection}}
    \renewcommand{\thesection}{\thechapter.\arabic{section}}
    \renewcommand{\theHsection}{\theHchapter.\arabic{section}}
  }
  {
    \setcounter{SaveTextSection}{\value{section}}
    \setboolean{appsection}{true}
    \ifthenelse{\equal{\thechapter}{\CurrentChapter}}{}
    {
      \setcounter{SaveAppSection}{0}
      \edef\CurrentChapter{\thechapter}
    }
    \setcounter{section}{\value{SaveAppSection}}
    \renewcommand{\thesection}{\thechapter.\Alph{section}}
    \renewcommand{\theHsection}{\theHchapter.\Alph{section}}
  }
}

\newcounter{SaveAppSubsection}
\newcounter{SaveTextSubsection}
\edef\CurrentSubsection{\thesubsection}

\newboolean{appsubsection}
\setboolean{appsubsection}{false}

\newcommand{\switchsubsection}
{
  \ifthenelse{
    \boolean{appsubsection}
  }
  {
    \setcounter{SaveAppSubsection}{\value{subsection}}
    \setboolean{appsubsection}{false}
    \ifthenelse{\equal{\thesection}{\CurrentSection}}{}
    {
      \setcounter{SaveTextSubsection}{0}
      \edef\CurrentSection{\thesection}
    }
    \setcounter{subsection}{\value{SaveTextSubsection}}
    \renewcommand{\thesubsection}{\thesection.\arabic{subsection}}
    \renewcommand{\theHsubsection}{\theHsection.\arabic{subsection}}
  }
  {
    \setcounter{SaveTextSubsection}{\value{subsection}}
    \setboolean{appsubsection}{true}
    \ifthenelse{\equal{\thesection}{\CurrentSection}}{}
    {
      \setcounter{SaveAppSubsection}{0}
      \edef\CurrentSection{\thesection}
    }
    \setcounter{subsection}{\value{SaveAppSubsection}}
    \renewcommand{\thesubsection}{\thesection.\Alph{subsection}}
    \renewcommand{\theHsubsection}{\theHsection.\Alph{subsection}}
  }
}

\newcounter{SaveAppSubsubsection}
\newcounter{SaveTextSubsubsection}

\newboolean{appsubsubsection}
\setboolean{appsubsubsection}{false}

\newcommand{\switchsubsubsection}
{
  \ifthenelse{
    \boolean{appsubsubsection}
  }
  {
    \setcounter{SaveAppSubsubsection}{\value{subsubsection}}
    \setboolean{appsubsubsection}{false}
    \ifthenelse{\equal{\thesubsection}{\CurrentSubsection}}{}
    {
      \setcounter{SaveTextSubsubsection}{0}
      \edef\CurrentSubsection{\thesubsection}
    }
    \setcounter{subsubsection}{\value{SaveTextSubsubsection}}
    \renewcommand{\thesubsubsection}{\thesubsection.\arabic{subsubsection}}
    \renewcommand{\theHsubsubsection}{\theHsubsection.\arabic{subsubsection}}
  }
  {
    \setcounter{SaveTextSubsubsection}{\value{subsubsection}}
    \setboolean{appsubsubsection}{true}
    \ifthenelse{\equal{\thesubsection}{\CurrentSubsection}}{}
    {
      \setcounter{SaveAppSubsubsection}{0}
      \edef\CurrentSubsection{\thesubsection}
    }
    \setcounter{subsubsection}{\value{SaveAppSubsubsection}}
    \renewcommand{\thesubsubsection}{\thesubsection.\Alph{subsubsection}}
    \renewcommand{\theHsubsubsection}{\theHsubsection.\Alph{subsubsection}}
  }
}

\makeatletter
\@ifclassloaded{book}{
\numberwithin{section}{chapter}
}{}
\makeatother
\numberwithin{figure}{section}
\numberwithin{table}{section}
\numberwithin{equation}{section}

\newtheorem{theorem}{Theorem}[section]
\newtheorem{lemma}[theorem]{Lemma}
\newtheorem{proposition}[theorem]{Proposition}
\newtheorem{corollary}[theorem]{Corollary}

\newtheorem{thm-def}[theorem]{Theorem-Definition}

\theoremstyle{definition}
\newtheorem{definition}[theorem]{Definition}
\newtheorem{example}[theorem]{Example}
\newtheorem{remark}[theorem]{Remark}

\theoremstyle{remark}

\newcommand{\Cat}{\mathrm{Cat}}

\newcommand{\CPM}{\mathrm{CPM}} % combinatorial premodel categories

\newcommand{\Set}{\mathrm{Set}}

\newcommand{\sSet}{\mathrm{sSet}}

\newcommand{\Top}{\mathrm{Top}}

\newdir{ >}{{}*!/-10pt/@{>}}

\makeatletter
\newcommand{\xRightarrow}[2][]{\ext@arrow 0359\Rightarrowfill@{\quad#1\quad}{#2}}
\newcommand{\xLeftarrow}[2][]{\ext@arrow 3228\Leftarrowfill@{\quad#1\quad}{#2}}
\makeatother

\DeclareMathOperator{\BL}{L} % Left Bousfield localization
\DeclareMathOperator{\BR}{R} % Right Bousfield localization
\DeclareMathOperator{\coker}{coker}

\DeclareMathOperator*{\colim}{colim}

\DeclareMathOperator{\Eq}{Eq}

\DeclareMathOperator{\Ho}{Ho}

\DeclareMathOperator{\map}{map}

\DeclareMathOperator{\mor}{mor}

\DeclareMathOperator{\PSh}{PSh}

\DeclareMathOperator{\Sh}{Sh}

\DeclareMathOperator{\sPSh}{sPSh}
\newcommand{\Un}[1]{\mathrm{U}(#1)} % Universal homotopy theory

\newcommand{\Vn}[1]{\mathrm{V}(#1)} % Co-Universal homotopy theory

\newcommand{\gangle}[1]{\left\langle{#1}\right\rangle}

\newcommand{\id}{\mathrm{id}}

\newcommand{\op}{\mathrm{op}}

\newcommand{\prfr}[2]{(#1)\(\,\Rightarrow\,\)(#2)}
\newcommand{\prfl}[2]{(#1)\(\,\Leftarrow\,\)(#2)}

\newcommand{\LI}{\mathrm{Liso}} 

\begin{document}

\title{Homotopical presentation of categories}

\author{Seunghun Lee}

\address{Department of Mathematics, Konkuk University,
Kwangjin-Gu Hwayang-dong 1, Seoul 143-701, Korea}

\email{mbrs@konkuk.ac.kr}

\thanks{This research was supported by Basic Science Research Program
  through the National Research Foundation of Korea(NRF) funded by the
  Ministry of Education(no. 2017R1D1A1B03027980)}

\subjclass[2020]{Primary 18N40; Secondary 14C20}

\date{\today}

\keywords{Universal model category, Bousfield localization,
  Quillen equivalence, Equivalence of category, Calculus of
  right fractions, Simplicial presheaf}

\begin{abstract}
  We give a criterion for a functor \(F:C\rightarrow B\)
  between small categories to generate a small presentation
  of the universal model category \(\Un B\) in the sense of
  Dugger.
\end{abstract}

\maketitle

\tableofcontents

\section{Introduction}
\label{sec:introduction-1}

The purpose of this note is to give a criterion for a
functor \(F:C\rightarrow B\) between small categories \(C\)
and \(B\) to generate a small presentation of the universal
model category \(\Un B\) associated with \(B\) in the sense
of Dugger (\cite{dugger-01}). We then use it to recover
three results well-known from the infinity categorical
perspective.

In \cite{lee-24-b}, we use it to give a model theoretic
characterization of the ampleness for a divisors on a
complex smooth projective varieties.

\begin{definition}
  \label{def:homotopical_presentation:1}
  Let \(F:C\rightarrow B\) be a functor between small
  categories \(C\) and \(B\). Let \(S\) be a set of
  morphisms in \(C\) mapped to an isomorphism in \(B\).  We
  say that the functor \(F\) \textbf{generates a homotopical
    presentation of the category \(B\) with a generator
    \(C\) and a relation \(S\)} if the Quillen adjunction
  \begin{equation}
    \label{eqn:presheaves_and_Q-divisors:135}
    F_{*}:\BL_{S}(\Un C)\rightleftarrows \Un B:F^{*}
  \end{equation}  
  associated with \(F\) is a Quillen equivalence where
%  \(\Un C\) and \(\Un B\) are the universal model categories
%  (\cite{dugger-01}) on \(C\) and \(B\) respectively and
  \(\BL_{S}(\Un C)\) is the left Bousfield localization
  (\cite{hirschhorn-03}) of \(\Un C\) at \(S\). When the set
  \(S\) in \eqref{eqn:presheaves_and_Q-divisors:135} is the
  set of all morphisms in \(C\) mapped to an isomorphism in
  \(B\), we simply say that the functor \(F\)
  \textbf{generates a homotopical presentation of the
    category \(B\)}.
\end{definition}

\begin{remark}
  \label{rem:homotopical_presentation:3}
  The Quillen equivalence
  \eqref{eqn:presheaves_and_Q-divisors:135} is what Dugger
  called a small presentation of the model category
  \(\Un B\).  The category \(\Un B\) is a homotopy colimit
  completion of \(B\), and the Quillen equivalence
  \eqref{eqn:presheaves_and_Q-divisors:135} presents
  \(\Un B\) with \(C\) as a category of generators and \(S\)
  as a set of relations.\footnote{As the author explains in
    Section~4 in \cite{dugger-01}, \(\Un B\) should have
    been called a co-universal model category because it is
    relevant to the study of homotopy colimits. If we adopt
    this term for \(\Un B\) then
    \eqref{eqn:presheaves_and_Q-divisors:135} should be
    called a co-homotopical presentation.}
\end{remark}

% The identity functor on \(B\) generates a trivial
% homotopical presentation of \(B\).  Below, we will provide
% three non-trivial examples, two in category theory and one
% in algebraic geometry.

Let \(\sSet=\Set^{\Delta^{\op}}\) be the category of
simplicial sets.  Given a small category \(A\), we denote by
\(\sPSh(A)\) the category of simplicial presheaves on \(A\).
The universal model category \(\Un A\) on \(A\) is the model
category of \(\sPSh(A)\) and the Bousfield-Kan model
structure; a fibration is an object-wise fibration in
\(\sSet\) and a weak equivalence is an object-wise weak
equivalence in \(\sSet\).  Given a functor
\(F:C\rightarrow B\) between small categories, there is an
adjunction
\begin{equation}
  \label{eqn:presheaves_and_Q-divisors:168}
  F_{*}:\sPSh(C)\rightleftarrows \sPSh(B):F^{*}
\end{equation}
obtained from \(F\).\footnote{We follow the notations in
  \cite{artin-62}. Instead of
  \eqref{eqn:presheaves_and_Q-divisors:168},
  \(F^{\dag}:\sPSh(C)\rightleftarrows \sPSh(B):F_{*}\) is used in
  Definition~2.3.1 in \cite{kashiwara-shapira-06}. We do not need the
  right Kan extension functor \(F^{\ddag}\).}  It is the underlying
adjunction of the Quillen adjunction
\begin{equation}
  \label{eqn:presheaves_and_Q-divisors:136}
  F_{*}:\Un C\rightleftarrows \Un B:F^{*}.
\end{equation}
The Quillen adjunction
\eqref{eqn:presheaves_and_Q-divisors:135} is obtained from
\eqref{eqn:presheaves_and_Q-divisors:136} by localizing
\(\Un C\) at \(S\).  Even if we localize with respect to all
morphisms in \(C\) mapped to an isomorphism in \(B\), we
need to impose a connectivity of fibers of \(F\) for
\eqref{eqn:presheaves_and_Q-divisors:135} to be a Quillen
equivalence. For this, we introduce the following notions.

\begin{definition}
  \label{def:homotopical_presentation:2}
  Let \(C,B\) be categories. Let \(F:C\rightarrow B\) be a
  functor.  Let \(d\) be an object of \(C\).
  \begin{enumerate}
  \item We denote by \(C_{d}\) the full subcategory of the
    comma category \((C\downarrow d)\) such that for each
    object \(\gangle{e,\phi:e\rightarrow d}\) of
    \((C\downarrow d)\), \(\gangle{e,\phi}\) is an object of
    \(C_{d}\) iff \(F\phi\) is an isomorphism in \(B\).
  \item For a subcategory \(A_{d}\) of \(C_{d}\), we say
    that \(F\) \textbf{has \(A_{d}\)-lifting property}
    (Definition~\ref{def:presheaves_and_Q-divisors:2}) if
    \begin{itemize}
    \item for each \(e\in C\) and a morphism
      \(\delta:Fd\rightarrow Fe\), there is an object
      \(\gangle{e',\phi:e'\rightarrow d}\) of \(A_{d}\) and
      a morphisms \(\psi:e'\rightarrow e\) in \(C\) such
      that \(\delta=F\psi \cdot (F\phi)^{-1}\) holds.
    \end{itemize}
  \item For a subcategory \(A_{d}\) of \(C_{d}\), we say
    that \(A_{d}\) is \textbf{\(F\)-admissible}
    (Definition~\ref{def:presheaves_and_Q-divisors:5}) if
    \(A_{d}\) is cofiltered and \(F\) has \(A_{d}\)-lifting
    property.
  \item For a family of subcategories \(A_{d}\) of \(C_{d}\)
    for \(d\in C\), we denote by \(S_{A}\) the set of
    morphisms \(\phi:e\rightarrow d\) in \(C\) such that
    \(\gangle{e,\phi}\in A_{d}\)
    (Definition~\ref{def:presheaves_and_Q-divisors:6}).
  \end{enumerate}
\end{definition}

The following is our main result.

\begin{theorem}
  \label{thm:presheaves_and_Q-divisors:1}
  Let \(B,C\) be small categories. Let \(F:C\rightarrow B\) be a
  functor. We assume that \(C\) has finite limits and \(F\) preserves
  finite limits.  We assume that the following hold.
  \begin{enumerate}
  \item \(F\) is essentially surjective on objects. 
  \item For each \(d\in C\), there is a \(F\)-admissible
    subcategory \(A_{d}\) of \(C_{d}\).
  \end{enumerate}
  Then \(F\) generates a homotopical presentation of the
  category \(B\) with a generator \(C\) and a relation
  \(S_{A}\).
  % I.e., the Quillen adjunction
  % \begin{equation}
  %   \label{eqn:presheaves_and_Q-divisors:48}
  %   F_{*}:\BL_{S_{A}}(\Un C)\rightleftarrows \Un B:F^{*}
  % \end{equation}
  % obtained from \eqref{eqn:presheaves_and_Q-divisors:136}
  % is a Quillen equivalence.
\end{theorem}

\begin{remark}
  \label{rem:presheaves_and_Q-divisors:46}
  In Theorem~\ref{thm:presheaves_and_Q-divisors:1}, we do
  not assume that \(A_{d}\) is full or \(A_{d}\) is
  functorial in \(d\).
\end{remark}

Below, we give three applications of
Theorem~\ref{thm:presheaves_and_Q-divisors:1}. They seem
well-known from the theory of \((\infty,1)\)-category theory
(\cite{cisinski-19},
\cite{Carranza-Kapulkin-Lindsey-23}). However, our proof is
direct without relying on \((\infty,1)\)-category theory.

It is known that \(C_{d}\) is cofiltered for each
\(d\in C\). The first application is the following
characterization.

\begin{corollary}
  \label{cor:presheaves_and_Q-divisors:1}
  Let \(B,C\) be small categories. Let \(F:C\rightarrow B\) be a
  functor. We assume that \(C\) has finite limits and \(F\) preserves
  finite limits.  Then the following are equivalent.
  \begin{enumerate}
  \item\label{item:20} The following two properties hold.
    \begin{itemize}
    \item \(F\) is essentially surjective on objects.
    \item For each \(d\in C\), \(F\) has \(C_{d}\)-lifting
      property.
    \end{itemize}
  \item\label{item:21} \(F\) generates a homotopical
    presentation of \(B\).
    % \textcolor{red}{Out: I.e., if we let \(S_{C}\) be the
    % set of all morphisms in \(C\) mapped to an isomorphsm
    % in \(B\) then the Quillen adjunction
    %\begin{equation}
    %  \label{eqn:homotopical_presentation:1}
    %  F_{*}:\BL_{S_{C}}(\Un C)\rightleftarrows \Un B:F^{*}
    %\end{equation}
    %obtained from \eqref{eqn:presheaves_and_Q-divisors:136}
    % is a Quillen equivalence.}
  \end{enumerate}
\end{corollary}

\begin{remark}
  \label{rem:presheaves_and_Q-divisors:29}
  In Corollary~\ref{cor:presheaves_and_Q-divisors:1}, if we
  let \(S_{C}\) be the set of all morphisms in \(C\) mapped
  to an isomorphism in \(B\) then \(S_{C}\) admits the
  calculus of right fractions in the sense of Gabriel-Zisman
  (\cite{gabriel-zisman-67})\footnote{There is a notion of
    "calculus of fractions" for higher categories introduced
    by Cisinski in~\cite{cisinski-19}. In this note, we
    always work with the calculus of fractions introduced by
    Gabriel-Zisman in~\cite{gabriel-zisman-67}.}. Thus \(F\)
  has \(C_{d}\)-lifting property for all \(d\in C\) iff the
  associated functor \(C[(S_{C})^{-1}]\rightarrow B\) is
  full.
\end{remark}

\begin{definition}
  \label{def:homotopical_presentation:4}
  We call \eqref{eqn:presheaves_and_Q-divisors:136} the
  \textbf{universal adjunction} of \(F\).  We call \(F\) a
  \textbf{universal equivalence} if the universal adjunction
  of \(F\) is a Quillen equivalence.
\end{definition}

The second application is a characterization of the
equivalences of categories in terms of universal
equivalences.  A universal equivalence is full and
faithful. Every full functor has \(C_{d}\)-lifting property
for all \(d\in C\). Every full and faithful functor
\(F:C\rightarrow B\) satisfies \(\BL_{S_{C}}\Un C=\Un C\).
So the following
Corollary~\ref{cor:presheaves_and_Q-divisors:6} fits into
the framework of the above theorem.

\begin{corollary}
  \label{cor:presheaves_and_Q-divisors:6}
  Let \(B,C\) be small categories.  Let \(F:C\rightarrow B\) be a
  functor. We assume that \(C\) has finite limits and \(F\) preserves
  finite limits. Then the following are equivalent.\footnote{We assume
    the axiom of choice.}
  \begin{enumerate}
  \item\label{item:18} \(F\) is an equivalence of categories.
  \item\label{item:19} \(F\) is a universal equivalence.
  \end{enumerate}
\end{corollary}

\begin{remark}
  \label{rem:presheaves_and_Q-divisors:19}
  Corollary~\ref{cor:presheaves_and_Q-divisors:6}, hence
  from Corollary~\ref{cor:presheaves_and_Q-divisors:1}, does
  not hold without a restriction on \(F\). For example, if
  \(F:C\rightarrow B\) is a Cauchy completion
  (\cite{borceux-dejean-86}) of a category \(C\), then the
  induced adjunction is an Quillen equivalence. If a
  category has equalizers then it is Cauchy complete.
  Placing a restriction with finite limits provides one
  sufficient condition. But, this example also shows that
  Corollary~\ref{cor:presheaves_and_Q-divisors:6} is not
  optimal.
\end{remark}

%\begin{remark}
%  \label{rem:presheaves_and_Q-divisors:33}
%  \sout{Corollary~\ref{cor:presheaves_and_Q-divisors:6} shows that
%  the \(C_{d}\)-lifting property in
%  Definition~\ref{def:homotopical_presentation:2}
%  generalizes the fully faithfulness.}
%\end{remark}

\begin{remark}
  \label{rem:presheaves_and_Q-divisors:28}
  There is a unique model structure on the category \(\Cat\)
  of small categories whose set of weak equivalences are
  precisely the set of equivalences of small categories
  (\cite{rezk-96},\cite{joyal-d},\cite{Mathoeverflow-2010}).
  And Quillen equivalences are to be the weak equivalences
  among model categories (\cite{hovey-99}).
  Corollary~\ref{cor:presheaves_and_Q-divisors:6} says that
  \(\Un -\) as a functor on \(\Cat\) reflects the
  homotopical aspect of small categories and equivalences of
  small categories. This statement can be made more precise
  using the combinatorial \(\sSet\)-premodel categories
  \(\CPM_{\sSet}\) (\cite{barton-19}).  The (\(2\)-)category
  of model categories lacks limits and colimits. So, Barton
  introduced premodel categories to rectify this problem. A
  premodel category is a model category without the two out
  of three property for the set of weak equivalences
  (cf.~\cite{henry-20},\cite{lee-15}).  Barton showed that
  \(\CPM_{\sSet}\) has a model \(2\)-category structure.
  Corollary~\ref{cor:presheaves_and_Q-divisors:6} says that
  the functor \(\Un -:\Cat\rightarrow \CPM_{\sSet}\)
  preserves and reflects the weak equivalences when
  restricted to the category of small categories with finite
  limits and finite limit preserving functors.
\end{remark}

The last application explains how \(\Un -\) acts on the
calculus of right fractions (\cite{gabriel-zisman-67}).  Let
\(C\) be a small category.  Let \(\Sigma\) be a set of
morphisms in \( C\) satisfying the calculus of right
fractions (\cite{gabriel-zisman-67}).  Then the localization
\(F:C\rightarrow C[\Sigma^{-1}]\) of \(C\) at \(\Sigma\) has
the \(C_{d}\)-lifting property for all \(d\in C\). So the
following corollary also fits into the framework of
Corollary~\ref{cor:presheaves_and_Q-divisors:1}.

\begin{corollary}
  \label{cor:presheaves_and_Q-divisors:2}
  Let \(C\) be a small category.  Let \(\Sigma\) be a subset of
  \(\mor C\) admitting the calculus of right fractions.  Let
  \(F:C\rightarrow C[\Sigma^{-1}]\) be the localization of \(C\) at
  \(\Sigma\). We assume that \(C\) has finite limits. Then the
  adjunction
  \begin{equation}
    \label{eqn:presheaves_and_Q-divisors:106}
    F_{*}:\BL_{\Sigma}(\Un C)\rightleftarrows \Un {C[\Sigma^{-1}]}:F^{*}
  \end{equation}
  is a Quillen equivalence.
\end{corollary}

\begin{remark}
  \label{rem:presheaves_and_Q-divisors:24}
  In Corollary~\ref{cor:presheaves_and_Q-divisors:2}, we are not
  assuming that \(\Sigma\) is saturated.
\end{remark}

In Section~\ref{sec:admiss-subc}, we prove a key property of
the admissible subcategory while fixing some notations. In
Section~\ref{sec:proof-theorem-ref} we prove the results
stated in this section. In Section~\ref{sec:examples}, we
give some examples of homotopical presentations of
categories.

% \subsection*{Acknowledgment}
% \label{sec:acknowledgment}

% In an attempt to publish a previous version of this paper,
% the anonymous referee pointed out that the proof of
% Corollary~\ref{cor:presheaves_and_Q-divisors:6} has a gap,
% and explained the example in
% Remark~\ref{rem:presheaves_and_Q-divisors:19}.  Referee's
% remark lead me to the current proof of the Lemma. So, I
% would like to thank the referee for the report.

\section{Review on Model category}
\label{sec:soutr-site-model}

Here we briefly review the model category, the universal
homotopy theory and the Bousfield localization, referring
the details to \cite{quillen-67}, \cite{dugger-01} and
\cite{hirschhorn-03}.

\subsection{Model categories}
\label{sec:model-categories}

A \textbf{model category}
(\cite{quillen-67},\cite{may-ponto-12})
consists of a category \(M\) with small limits and small colimits and
three sets \(W,C,F\) of morphisms such that
\begin{enumerate}
\item \(W\) satisfies the two out of three properties.
\item \((C\cap W,F)\) and \((C,F\cap W)\) are weak factorization systems.
\end{enumerate}
Morphisms in \(W\), \(C\) and \(F\) are called \textbf{weak
  equivalences, cofibrations and fibrations} respectively. A morphism
is called a \textbf{trivial fibration} if it is a weak equivalence and
a fibration.  A morphism is called a \textbf{trivial cofibration} if
it is a weak equivalence and a cofibration.  An object is
\textbf{fibrant} if the unique morphism to the terminal object is a
fibration.  An object is \textbf{cofibrant} if the unique morphism from
the initial object is a cofibration.  Given a object \(x\) of a model
category we denote by
\begin{equation}
  \label{eqn:presheaves_and_Q-divisors:77}
  p_{x}:Qx\rightarrow x
\end{equation}
a trivial fibration from a cofibrant object  \(Qx\) of \(M\), and by
\begin{equation}
  \label{eqn:presheaves_and_Q-divisors:78}
  i_{x}:x\rightarrow Rx
\end{equation}
a trivial cofibration from a fibrant object \(Rx\) of \(M\).  We call
\eqref{eqn:presheaves_and_Q-divisors:77} a \textbf{cofibrant
  replacement} of \(x\) and \eqref{eqn:presheaves_and_Q-divisors:78} a
\textbf{fibrant replacement} of \(x\).

\subsection{Universal homotopy theory}
\label{sec:univ-homot-theory}

Let \(\sSet\) be the category of simplicial sets.
Let \(C\) be a small category.
The category of the presheaves of simplicial sets on \(C\)
\begin{equation}
  \label{eqn:presheaves_and_Q-divisors:21}
  \sPSh(C)=\sSet^{C^{\op}}
\end{equation}
has a proper simplicial model structure, called the
\textbf{Bousfield-Kan model structure}, where the weak equivalences
and the fibrations are defined object-wise (p.~314 in
\cite{bousfield-kan}).  Following \cite{dugger-01}, we denote by
\begin{equation}
  \label{eqn:presheaves_and_Q-divisors:128}
  \Un C
\end{equation}
the model category of the category \(\sPSh(C)\) of simplicial
presheaves on \(C\) equipped with the Bousfield-Kan model structure.

Let \(F:C\rightarrow B\) a functor between small categories.  There is
an adjunction (Chapter I.2 in \cite{artin-62}, Chapter~2.3
in~\cite{kashiwara-shapira-06})\footnote{\(F^{*}\) is \(F_{*}\) in
  \cite{kashiwara-shapira-06}.}
\begin{equation}
  \label{eqn:presheaves_and_Q-divisors:16}
  F_{*}:\sPSh(C)\rightleftarrows \sPSh(B):F^{*}
\end{equation}
where, for each \(Y\in \sPSh(B)\),
\begin{equation}
  \label{eqn:presheaves_and_Q-divisors:17}
  F^{*}Y=Y\cdot F
\end{equation}
and, for each \(X\in \sPSh(C)\) and \(b\in B\),
\begin{equation}
  \label{eqn:presheaves_and_Q-divisors:18}
  F_{*}X(b)=\colim_{\substack{\gangle{c,b\rightarrow Fc}\\\in (b\downarrow F)^{\op}}} X(c)
\end{equation}
where \((b\downarrow F)\) is the comma category.\footnote{The right
  hand side is the colimit of
  \begin{equation}
    \label{eqn:presheaves_and_Q-divisors:129}
    X\cdot P : (b\downarrow F)^{\op}\rightarrow \sSet
  \end{equation}
  where \(P:(b\downarrow F)^{\op}\rightarrow C\) mapping
  \(\gangle{c,b \rightarrow Fc}\)}
The adjunction
\eqref{eqn:presheaves_and_Q-divisors:16} is the Quillen adjunction
%\begin{equation}
%  \label{eqn:presheaves_and_Q-divisors:31}
%  \Re:\Un C\rightleftarrows \Un B:\Sing
%\end{equation}
associated with the cosimplicial resolution
\begin{equation}
  \label{eqn:presheaves_and_Q-divisors:25}
  \Gamma:C\times\Delta\rightarrow \Un B, \qquad
  \Gamma(c,[n])=r_{B}(Fc)\otimes \Delta[n]
\end{equation}
of \(r_{B}\cdot F\) (\cite{dugger-01}).  The
adjunction~\eqref{eqn:presheaves_and_Q-divisors:16} becomes
a Quillen adjunction between \(\Un C\) and \(\Un B\).
\begin{equation}
  \label{eqn:presheaves_and_Q-divisors:167}
  F_{*}:\Un C\rightleftarrows \Un B:F^{*}
\end{equation}

Given a category \(C\), we denote by \(r_{C}\) the functor
\begin{equation}
  \label{eqn:presheaves_and_Q-divisors:104}
  r_{C}:C\rightarrow\sPSh(C)
\end{equation}
mapping an object \(d\) of \(C\) to the simplicial presheaf on \(C\)
whose \(n\)-th component is  \(\hom_{C}(-,d)\) for all \(n\in\Delta\)
and the face and the degeneracy maps are all identities.
By an abuse
of notation, we also denote by \(r_{C}\) the Yoneda embedding
\begin{equation}
  \label{eqn:presheaves_and_Q-divisors:284}
  r_{C}:C\rightarrow \PSh(C).
\end{equation}
We also denote \(r_{C}\) by \(r\) for simplicity.

%Let
%\begin{equation}
%  \label{eqn:presheaves_and_Q-divisors:22}
%  r_{C}:C\rightarrow \Un C
%\end{equation}
%be the functor that maps an object \(c\) of \(C\) to the constant
%simplicial presheaf \(\hom_{C}(-,c)\). We also denote \(r_{C}\) by
%\(r\) for simplicity.  

We need the following well-known lemma. It does not requires
finite limits. See, for example, (1.1.10.3) in
Theorem~1.1.10 in~\cite{cisinski-19}.

\begin{lemma}[Theorem~I.1.10 in~\cite{cisinski-19}]
  \label{lem:presheaves_and_Q-divisors:20}
  Let \(F:C\rightarrow B\) a functor between locally small categories.
  Let \(y\) be an object of \(C\).
  Then, there is an canonical isomorphism 
  \begin{equation}
    \label{eqn:presheaves_and_Q-divisors:139}
    F_{*}\hom_{C}(-,y)\xrightarrow[\cong]{\alpha} \hom_{B}(-,Fy).
  \end{equation}
\end{lemma}

\subsection{Bousfield localization}
\label{sec:bousf-local}

Let \(M\) be a model category. Let \(S\) be a set of morphisms in
\(M\). An object \(Z\) of \(M\) is called a \(S\)-local object of
\(M\) if
\begin{enumerate}
\item \(Z\) is a fibrant object of \(M\).
\item For every morphism \(f:X\rightarrow Y\) in \(S\), the induced
  morphism between function complexes (Ch.~17 in \cite{hirschhorn-03})
  \begin{equation}
    \label{eqn:presheaves_and_Q-divisors:34}
    f^{*}:\map(Y,Z)\rightarrow \map(X,Z)
  \end{equation}
  is a weak equivalence in \(\sSet\).
\end{enumerate}
A morphism \(f:X\rightarrow Y\)  in \(M\) is a \(S\)-local equivalence
if for every \(S\)-local object \(Z\) of \(M\), the induced morphism
\begin{equation}
  \label{eqn:presheaves_and_Q-divisors:35}
  f^{*}:\map(Y,Z)\rightarrow \map(X,Z)
\end{equation}
is a weak equivalence in \(\sSet\).  The \textbf{left Bousfield
  localization} \(\BL_{S}(M)\) of \(M\) with respect to \(S\) is a model
category such that
\begin{enumerate}
\item The underlying category is \(M\).
\item The set of weak equivalences is the set of \(S\)-local
  equivalences in \(M\).
\item The set of cofibrations is the set of cofibrations in \(M\).
\end{enumerate}

For a small category \(C\) and a set \(S\) of morphisms in \(C\) the
following holds (\cite{hirschhorn-03}).  By an abuse of notation, we
also denote by \(S\) its image \(r_{C}(S)\).
\begin{enumerate}
\item The left Bousfield localization \(\BL_{S}\Un C\) exists.
\item \(\BL_{S}\Un C\) inherits a left proper simplicial model structure
  from \(\Un C\).
\item The fibrant objects are precisely the \(S\)-local objects in
  \(\Un C\).
\end{enumerate}
Because \(\Un C\) is a simplicial model category, if \(X\) and \(Y\) are
cofibrant objects then we can use the enriched hom-sets
\(\underline{\hom}(X,Z)\) and \(\underline{\hom}(Y,Z)\) instead of the
function complexes \(\map(X,Z)\) and \(\map(Y,Z)\) respectively to
define \(S\)-local objects.  So by Yoneda lemma, we obtain the
following lemma.

\begin{lemma}
  \label{lem:presheaves_and_Q-divisors:6}
  Let \(C\) be a small category. Let \(S\) be a set of morphisms in
  \(C\). Then, for every fibrant object \(X\) of \(\Un C\), the following
  are equivalent.
  \begin{enumerate}
  \item \(X\) is a \(S\)-local object of \(\Un C\).
  \item For every morphism \(f:c_{1}\rightarrow c_{2}\) in \(S\), the
    induced morphism
    \begin{equation}
      \label{eqn:presheaves_and_Q-divisors:36}
      f^{*}:X(c_{2})\rightarrow X(c_{1})
    \end{equation}
    is a weak equivalence in \(\sSet\).
  \end{enumerate}
\end{lemma}

\section{Admissible subcategory}
\label{sec:admiss-subc}

In this section, we prove a key property
Lemma~\ref{lem:presheaves_and_Q-divisors:2} of admissible
subcategories
(Definition~\ref{def:homotopical_presentation:2}). We restate
Definition~\ref{def:homotopical_presentation:2} for the
convenience of readers.  And we set some notations that will
be used throughout this note.

The lifting property in
Definition~\ref{def:homotopical_presentation:2} is related
with the initial functor and the calculus of right
fractions. They are explained in the two following sections.

\subsection{Initial functor and admissible subcategory}
\label{sec:init-funct-admiss}

We begin with the key notions.

\begin{definition}
  \label{def:presheaves_and_Q-divisors:1}
  Let \(B,C\) be categories.
  Let \(F:C\rightarrow B\) be a functor.
  Let \(d\) be an object of \(C\).
  \begin{enumerate}
  \item We denote by \(C_{d}\) the full subcategory of the comma
    category \((C\downarrow d)\) such that for each object
    \(\gangle{e,\phi:e\rightarrow d}\) of \((C\downarrow d)\),
    \(\gangle{e,\phi}\) is an object of \(C_{d}\) iff \(F\phi\) is an
    isomorphism in \(B\).
  \item We denote by \(B_{d}\) the comma category
    \((Fd\downarrow F)\). Thus, the objects are the pairs
    \begin{equation}
      \label{eqn:presheaves_and_Q-divisors:11}
      \gangle{e,\delta:Fd\rightarrow Fe}
    \end{equation}
    of an object \(e\) of \(C\) and a morphism
    \(\delta:Fd\rightarrow Fe\) in \(B\), and the morphisms
    \(\psi:\gangle{e_{1},\delta_{1}}\rightarrow
    \gangle{e_{2},\delta_{2}}\) are the morphisms
    \(\psi:e_{1}\rightarrow e_{2}\) in \(C\) such that
    \(F\psi\cdot\delta_{1}=\delta_{2}\), i.e. the diagram in \(B\)
    \begin{equation}
      \label{eqn:presheaves_and_Q-divisors:10}
        \begin{tikzcd}
          Fe_{1} \ar[rr,"F\psi"] && Fe_{2} \\
          & Fd \ar[ul,"\delta_{1}"]
          \ar[ur,"\delta_{2}"'] &
        \end{tikzcd}
    \end{equation}
    commutes.
  \item We denote by \(F_{d}:C_{d}\rightarrow B_{d}\) the functor such
    that
    \begin{enumerate}
    \item \(F_{d}\gangle{e,\phi}=\gangle{e,(F\phi)^{-1}}\) for each
      object \(\gangle{e,\phi}\) of \(C_{d}\).
    \item \(F_{d}\psi=\psi\) for each morphism
      \(\gangle{e_{1},\phi_{1}}\xrightarrow\psi\gangle{e_{2},\phi_{2}}\)
      in \(C_{d}\).
    \end{enumerate}
    Thus, the functor \(F_{d}\) maps a commutative diagram
    \begin{equation}
      \label{eqn:presheaves_and_Q-divisors:1}
      \begin{tikzcd}
        e_{1} \ar[rr,"\psi"] \ar[dr,"\phi_{1}"'] && e_{2} \ar[dl,"\phi_{2}"]\\
        & d &
      \end{tikzcd}
    \end{equation}
    in \(C\) to a commutative diagram
    \begin{equation}
      \label{eqn:presheaves_and_Q-divisors:15}
      \begin{tikzcd}
        Fe_{1} \ar[rr,"F\psi"] && Fe_{2} \\
        & Fd \ar[ul,"(F\phi_{1})^{-1}","\cong"']
        \ar[ur,"(F\phi_{2})^{-1}"', "\cong"] &
      \end{tikzcd}
    \end{equation}
    in \(B\).
  \end{enumerate}
\end{definition}

\begin{remark}
  \label{rem:presheaves_and_Q-divisors:32}
  \(C_{d},B_{d}\) depend on \(F\) and should be written as
  \(C_{d,F},B_{d,F}\). But, throughout this note, we always
  work with a fixed functor \(F\). So we will omit \(F\) for
  simplicity. 
\end{remark}

\begin{remark}
  \label{rem:presheaves_and_Q-divisors:1}
  If \(\psi\) is a morphism in \(C_{d}\) then \(F\psi\) is an
  isomorphism in \(B\). Thus the functor \(F_{d}\) is well-defined.
\end{remark}

We recall the (co)filtered categories.

\begin{definition}
  \label{def:presheaves_and_Q-divisors:3}
  A non-empty category \(J\) is called \textbf{cofiltered} if the
  following hold.
  \begin{enumerate}
  \item\label{item:22} To any two objects \(j,j'\in J\), there is
    \(i\in J\) and morphisms \(i\rightarrow j\) and
    \(i\rightarrow j'\) in \(J\).
  \item\label{item:23} To any two parallel morphisms
    \(u,v:j\rightrightarrows k\), there is \(i\in J\) and morphism
    \(w:i\rightarrow j\) such that \(uw=vw\).
  \end{enumerate}
 The category \(J\) is called \textbf{filtered} if \(J^{\op}\) is cofiltered.  
\end{definition}

The following lemma is well-known.

\begin{lemma}
  \label{lem:presheaves_and_Q-divisors:30}
  Let \(B,C\) be categories. Let \(F:C\rightarrow B\) be a functor. We
  assume that \(C\) has finite limits and \(F\) preserves finite
  limits.  Then for each object \(b\) of \(B\), the comma category
  \((b\downarrow F)\) is cofiltered.
\end{lemma}

\begin{lemma}
  \label{lem:presheaves_and_Q-divisors:1}
  Let \(B,C\) be categories. Let \(F:C\rightarrow B\) be a functor. We
  assume that \(C\) has finite limits and \(F\) preserves finite
  limits.  Then for each object \(d\) of \(C\), the category \(C_{d}\)
  is cofiltered.
\end{lemma}
\begin{proof}
  Let \(d\) be an object of \(C\).

  Let \(\gangle{e_{1},\phi_{1}}\) and \(\gangle{e_{2},\phi_{2}}\) be
  two objects of \(C_{d}\). We assume that the left of
  \eqref{eqn:presheaves_and_Q-divisors:4} is the pullback
  diagram. Since \(F\) preserves finite limits, the diagram on the
  right is also the pullback diagram.
  \begin{equation}
    \label{eqn:presheaves_and_Q-divisors:4}
    \begin{tikzcd}
      e \ar[r,"\pi_{2}"] \ar[d,"\pi_{1}"'] & e_{2} \ar[d,"\phi_{2}"]\\
      e_{1} \ar[r,"\phi_{1}"] & d
    \end{tikzcd}
    \qquad
    \begin{tikzcd}
      Fe \ar[r,"F\pi_{2}"] \ar[d,"F\pi_{1}"'] & Fe_{2} \ar[d,"F\phi_{2}"]\\
      Fe_{1} \ar[r,"F\phi_{1}"] & Fd
    \end{tikzcd}    
  \end{equation}
  Since \(F\phi_{1}\) and \(F\phi_{2}\) are isomorphisms in \(B\),
  so are \(F\pi_{1}\) and \(F\pi_{2}\). Let
  \(\phi=\phi_{1}\cdot\pi_{1}\). Then
  \(\gangle{e,\phi}\) is an object of \(C_{d}\) and we have morphisms
  \begin{equation}
    \label{eqn:presheaves_and_Q-divisors:5}
    \pi_{i}:\gangle{e,\phi}\rightarrow\gangle{e_{i},\phi_{i}}
  \end{equation}
  in \(C_{d}\) for \(i=1,2\).

  Let
  \(\psi_{1},\psi_{2}:\gangle{e_{1},\phi_{1}} \rightrightarrows
  \gangle{e_{2},\phi_{2}}\) be two morphisms in \(C_{d}\). We assume
  that the left of \eqref{eqn:presheaves_and_Q-divisors:6} is the the
  equalizer of \(\psi_{1},\psi_{2}\). Since \(F\) preserves finite
  limits, the diagram on the right is also the equalizer of
  \(F\psi_{1},F\psi_{2}\).
  \begin{equation}
    \label{eqn:presheaves_and_Q-divisors:6}
    \begin{tikzcd}
      e\ar[r,"\varepsilon"] & e_{1}\ar[r, shift left, "\psi_{1}"]
      \ar[r, shift right, "\psi_{2}"'] & e_{2}
    \end{tikzcd}
    \qquad
    \begin{tikzcd}
      Fe\ar[r,"F\varepsilon"] & Fe_{1}\ar[r, shift left, "F\psi_{1}"]
      \ar[r, shift right, "F\psi_{2}"'] & Fe_{2}
    \end{tikzcd}
  \end{equation}
  Because \(F\phi_{2}\) is an isomorphism, \(F\varepsilon\) is also an
  isomorphism. Let \(\phi=\phi_{1}\cdot\varepsilon\). Then
  \(\gangle{e,\phi}\) is an object of \(C_{d}\) and
  \begin{equation}
    \label{eqn:presheaves_and_Q-divisors:7}
    \psi_{1}\cdot\varepsilon=\psi_{2}\cdot\varepsilon:
    \gangle{e,\phi}\rightarrow \gangle{e_{2},\psi_{2}}
  \end{equation}
  holds in \(C_{d}\).
\end{proof}

% The following property (Lift) naturally occurs in
% Section~\ref{sec:proofs-corollaries}.

The following definition of \(A_{d}\)-lifting property is
equivalent to that of
Definition~\ref{def:homotopical_presentation:2}.

\begin{definition}
  \label{def:presheaves_and_Q-divisors:2}
  Let \(B,C\) be categories.  Let \(F:C\rightarrow B\) be a
  functor.  Let \(d\) be an object of \(C\).  Let \(A_{d}\)
  be a non-empty subcategory of \(C_{d}\).  We say that
  \(F\) \textbf{have \(A_{d}\)-lifting property} if for each
  \(\gangle{e,\delta}\in B_{d}\), the comma category
  \((F_{d}|_{A_{d}}\downarrow \gangle{e,\delta})\) is not
  empty where \(F_{d}|_{A_{d}}:A_{d}\rightarrow B_{d}\) is
  the restriction.
  \begin{equation}
    \label{eqn:presheaves_and_Q-divisors:3}
    (F_{d}|_{A_{d}}\downarrow \gangle{e,\delta})\not=\emptyset
  \end{equation}
  In other words, for each \(\gangle{e,\delta}\in B_{d}\), there is an
  object \(e'\) of \(C\) and morphisms \(\phi:e'\rightarrow d\) and
  \(\psi:e'\rightarrow e\) in \(C\) such that
  \begin{enumerate}
  \item \(\gangle{e',\phi}\) is an object of \(A_{d}\).
  \item \(\psi:F_{d}\gangle{e',\phi}\rightarrow \gangle{e,\delta}\) is
    a morphism in \(B_{d}\), i.e. the diagram
    \begin{equation}
      \label{eqn:presheaves_and_Q-divisors:2}
      \begin{tikzcd}
        Fe' \ar[rr,"F\psi"] && Fe \\
        & Fd \ar[ul,"F(\phi)^{-1}","\cong"']
        \ar[ur,"\delta"'] &
      \end{tikzcd}
      \quad
      (\delta\cdot F(\phi)=F(\psi))
    \end{equation}
    in \(B\) commutes.
  \end{enumerate}
\end{definition}

\begin{definition}
  \label{def:presheaves_and_Q-divisors:5}
  Let \(B,C\) be categories.  Let \(F:C\rightarrow B\) be a
  functor.  Let \(d\) be an object of \(C\).  A non-empty
  subcategory \(A_{d}\) of \(C_{d}\) is called
  \textbf{\(F\)-admissible} if
  \begin{enumerate}
  \item \(A_{d}\) is cofiltered.
  \item \(F\) have \(A_{d}\)-lifting property.
  \end{enumerate}
\end{definition}

\begin{example}
  \label{exa:presheaves_and_Q-divisors:1}
  Let \(B,C\) be categories. Let \(F:C\rightarrow B\) be a
  functor. We assume that \(C\) has finite limits and \(F\)
  preserves finite limits.  Let \(d\) be an object of \(C\).
  Then, by Lemma~\ref{lem:presheaves_and_Q-divisors:1}, the
  following are equivalent.
  \begin{enumerate}
  \item \(C_{d}\) is \(F\)-admissible.
  \item \(F\) have \(C_{d}\)-lifting property.
  \end{enumerate}
\end{example}

We recall the initial functors (Chapter~IX.3 in~\cite{maclane-98}).

\begin{definition}
  \label{def:presheaves_and_Q-divisors:4}
  A functor \(L:J'\rightarrow J\) is called initial if for each
  \(k\in J\) the comma category \((L\downarrow k)\) is non-empty and
  connected.  \(L\) is called final if \(L^{\op}\) is initial.
\end{definition}

The following is the key property of \(F\)-admissible
subcategories.

\begin{lemma}
  \label{lem:presheaves_and_Q-divisors:2}
  Let \(B,C\) be categories. Let \(F:C\rightarrow B\) be a
  functor. We assume that \(C\) has finite limits and \(F\)
  preserves finite limits.  Let \(d\) be an object of \(C\).
  If \(A_{d}\) is a \(F\)-admissible then the restriction
  \(F_{d}|_{A_{d}}\) is an initial functor.
\end{lemma}
\begin{proof}
  Let \(\gangle{e,\delta}\) be an object of \(B_{d}\).  The comma
  category \((F_{d}|_{A_{d}}\downarrow \gangle{e,\delta})\) is
  non-empty by the assumption. So, we only need to show that
  \((F_{d}|_{A_{d}}\downarrow \gangle{e,\delta})\) is connected.

  Consider two objects \(\psi_{1},\psi_{2}\) of
  \((F_{d}|_{A_{d}}\downarrow \gangle{e,\delta})\).
  \begin{equation}
    \label{eqn:presheaves_and_Q-divisors:8}
    \begin{tikzcd}
      F_{d}\gangle{e_{1},\phi_{1}} \ar[dr,"\psi_{1}"'] &&
      F_{d}\gangle{e_{2},\phi_{2}} \ar[dl,"\psi_{2}"]\\
      & \gangle{e,\delta}
    \end{tikzcd}
  \end{equation}
  We will complete the diagram into a commutative diagram 
  \begin{equation}
    \label{eqn:presheaves_and_Q-divisors:9}
    \begin{tikzcd}
      F_{d}\gangle{e_{1},\phi_{1}} \ar[dr,"\psi_{1}"'] &
      F_{d}\gangle{e_{3},\phi_{3}} \ar[l,"\rho_{1}"', dotted]
      \ar[d,"\psi_{3}" ] \ar[r,"\rho_{2}", dotted]&
      F_{d}\gangle{e_{2},\phi_{2}} \ar[dl,"\psi_{2}"]\\
      & \gangle{e,\delta}
    \end{tikzcd}    
  \end{equation}
  in \((F_{d}|_{A_{d}}\downarrow \gangle{e,\delta})\) with another
  object
  \(\psi_{3}:F_{d}\gangle{e_{3},\phi_{3}}\rightarrow
  \gangle{e,\delta}\) of
  \((F_{d}|_{A_{d}}\downarrow \gangle{e,\delta})\) and two morphisms
  \(\rho_{1}:\psi_{3}\rightarrow \psi_{1}\) and
  \(\rho_{2}:\psi_{3}\rightarrow \psi_{2}\) in
  \((F_{d}|_{A_{d}}\downarrow \gangle{e,\delta})\), i.e.  two
  morphisms \(\phi_{3}:e_{3}\rightarrow d\) and
  \(\psi_{3}:e_{3}\rightarrow e\) in \(C\) satisfying
  \begin{equation}
    \label{eqn:presheaves_and_Q-divisors:14}
    \gangle{e_{3},\phi_{3}}\in A_{d}
  \end{equation}
  and 
  \begin{equation}
    \label{eqn:presheaves_and_Q-divisors:102}
    \delta\cdot (F\phi_{3})=F\psi_{3}
  \end{equation}
  in \(B\), and two morphisms \(\rho_{1}:e_{3}\rightarrow e_{1}\) and
  \(\rho_{2}:e_{3}\rightarrow e_{2}\) in \(C\) satisfying
  \begin{align}
    \label{eqn:presheaves_and_Q-divisors:100}
    F\phi_{1}\cdot F\rho_{1}&=F\phi_{3}=F\phi_{2}\cdot F\rho_{2}\\
    \label{eqn:presheaves_and_Q-divisors:108}
    \psi_{1}\cdot\rho_{1}&=\psi_{3}=\psi_{2}\cdot\rho_{2}.
  \end{align}
  Notice that \eqref{eqn:presheaves_and_Q-divisors:102} follows from
  \eqref{eqn:presheaves_and_Q-divisors:100} and
  \eqref{eqn:presheaves_and_Q-divisors:108} because of
  \(\delta\cdot (F\phi_{1})=F\psi_{1}\).

  Consider the following diagram in \(C\).  
  \begin{equation}
    \label{eqn:presheaves_and_Q-divisors:109}
    \begin{tikzcd}
      & d &\\
      e_{1} \ar[ur,"\phi_{1}"] \ar[dr,"\psi_{1}"'] &
      \ar[l,"\rho_{1}'"'] e' \ar[r,"\rho_{2}'"] \ar[u,"\phi'"]
      & e_{2} \ar[ul,"\phi_{2}"'] \ar[dl,"\psi_{2}"] \\
      & e &
    \end{tikzcd}
  \end{equation}
  The two upper triangles are commutative diagrams that the
  connectivity property of the cofilteredness of \(A_{d}\) provides us
  with, hence \(\gangle{e',\phi'}\in A_{d}\).  The lower triangle may
  not commute, but
  \begin{equation}
    \label{eqn:presheaves_and_Q-divisors:171}
    F\psi_{1}\cdot F\rho_{1}'=
    F\psi_{2}\cdot F\rho_{2}'
  \end{equation}
  holds because of \(F\psi_{1}=\delta\cdot F\phi_{1}\) and
  \(F\psi_{2}=\delta\cdot F\phi_{2}\).

  Let \(\varepsilon':e''\rightarrow e'\) be the equalizer of
  \(\psi_{1}\cdot\rho_{1}'\) and \(\psi_{2}\cdot\rho_{2}'\).
  Since \(F\) preserves finite limits, \(F\varepsilon'\) is
  an isomorphism in \(B\) by
  \eqref{eqn:presheaves_and_Q-divisors:171}. Because
  \(F\phi'\) is also an isomorphism in \(B\), by applying
  the \(A_{d}\)-lifting property of \(F\) to
  \((F(\phi'\cdot\varepsilon'))^{-1}\), there is a (not
  necessarily commutative) diagram in \(C\)
  \begin{equation}
    \label{eqn:presheaves_and_Q-divisors:110}
    \begin{tikzcd}
      & e' \ar[dr,"\phi'"] &\\
      e'' \ar[ur,"\varepsilon'"] && d\\
      & e_{3} \ar[ul,"\varepsilon''"] \ar[ur,"\phi_{3}"'] &
    \end{tikzcd}
  \end{equation}
  such that
  \begin{equation}
    \label{eqn:presheaves_and_Q-divisors:172}
    \gangle{e_{3},\phi_{3}}\in A_{d}
  \end{equation}
  and
  \begin{equation}
    \label{eqn:presheaves_and_Q-divisors:105}
    F\phi_{3}=F\phi'\cdot F\varepsilon'\cdot F\varepsilon''
  \end{equation}
  hold.  We let
  \begin{align}
    \label{eqn:presheaves_and_Q-divisors:170}
    \rho_{1}&=\rho_{1}'\cdot \varepsilon'\cdot\varepsilon''\\
    \rho_{2}&=\rho_{2}'\cdot \varepsilon'\cdot\varepsilon''\\
    \psi_{3}&=\psi_{1}\cdot\rho_{1}(=\psi_{2}\cdot\rho_{2}).
  \end{align}
  Then \eqref{eqn:presheaves_and_Q-divisors:100} holds because of
  \eqref{eqn:presheaves_and_Q-divisors:105}.
  \eqref{eqn:presheaves_and_Q-divisors:108} holds by definition.
\end{proof}

\begin{remark}
  \label{rem:presheaves_and_Q-divisors:14}
  In Lemma~\ref{lem:presheaves_and_Q-divisors:2}, if \(A_{d}=C_{d}\)
  for all \(d\in C\) then the converse also holds. Thus the
  following are equivalent.
  \begin{enumerate}
  \item \(F_{d}\) is initial.
  \item \(F\) has \(C_{d}\)-lifting property.
  \end{enumerate}
\end{remark}

\begin{definition}
  \label{def:presheaves_and_Q-divisors:6}
  Let \(C,B\) be categories. Let \(F:C\rightarrow B\) be a
  functor.
  \begin{enumerate}
  \item Let \(d\) be an object of \(C\). Let \(A_{d}\) be a
    subcategory of \(C_{d}\). We denote by \(S_{A_{d}}\) the
    set of all morphisms \(\phi:e\rightarrow d\) in \(C\)
    such that \(\gangle{e,\phi}\in A_{d}\).
    \begin{equation}
      \label{eqn:presheaves_and_Q-divisors:96}
      S_{A_{d}}=\bigl\{\phi:e\rightarrow d \text{ in } C \mid
      \gangle{e,\phi}\in A_{d}\bigr\}
    \end{equation}
  \item Assume that we are given a subcategory \(A_{d}\) for
    each \(d\). We denote by \(S_{A}\) the union of
    \(S_{A_{d}}\) over \(d\in C\).
    \begin{equation}
      \label{eqn:presheaves_and_Q-divisors:283}    
      S_{A}=\bigcup_{d\in C}S_{A_{d}}
    \end{equation}
    In particular, \(S_{C}\) is the set of all morphisms in
    \(C\) mapped to an isomorphism in \(B\).
    \begin{equation}
      \label{eqn:homotopical_presentation:3}
      S_{C}=\bigcup_{d\in C}S_{C_{d}}
    \end{equation}
  \end{enumerate}
\end{definition}

\subsection{Calculus of right fractions and admissible
  subcategory}
\label{sec:calc-fract-admiss}

We need the following observation in the proof of
Corollary~\ref{cor:presheaves_and_Q-divisors:1}.

\begin{lemma}
  \label{lem:homotopical_presentation:1}
  Let \(B,C\) be categories. Let \(F:C\rightarrow B\) be a
  functor. Let \(\Sigma\) be a set of morphisms in \(C\)
  such that
  \begin{enumerate}
  \item Every morphism in \(\Sigma\) mapped to an
    isomorphism in \(B\).
  \item \(\Sigma\) admits the calculus of right fractions.
  \end{enumerate}
  Let \(H:C[\Sigma^{-1}]\rightarrow B\) be the functor
  associated with \(F\). If \(H\) is full then \(F\) has
  \(C_{d}\)-lifting property for all \(d\in C\).
\end{lemma}

When \(\Sigma=S_{A}\) holds for a family of subcategories
\(A_{d}\) of \(C_{d}\), the following more precise result
holds.

\begin{lemma}
  \label{rem:presheaves_and_Q-divisors:47}
  Let \(C,B\) be categories. Let \(F:C\rightarrow B\) a
  functor. We assume that we are give a subcategory
  \(A_{d}\) of \(C_{d}\) for each object \(d\) of \(C\). Let
  \(H:C[S_{A}^{{-1}}]\rightarrow B\) be the associated functor.
  \begin{enumerate}
  \item If \(F\) has \(A_{d}\)-lifting property for all
    \(d\in C\) then \(H\) is full.
  \item If \(S_{A}\) admits the calculus of right fractions
    and \(H\) is full then \(F\) has \(A_{d}\)-lifting
    property for all \(d\in C\).
  \end{enumerate}
\end{lemma}
\begin{proof}
  (1) is clear.

  (2) Since \(S_{A}\) admits the calculus of right
  fractions, \((S_{A})^{-1}C=C[(S_{A})^{-1}]\). Thus \(F\)
  have \(A_{d}\)-lifting property for all \(d\in C\) by
  Proposition~I.2.4 in~\cite{gabriel-zisman-67}.
\end{proof}

The following Lemma~\ref{lem:homotopical_presentation:2}
shows that we can obtain an admissible subcategories
\(A^{\Sigma}_{d}\) from a set \(\Sigma\) of morphisms
admitting the calculus of right fractions.

\begin{definition}
  \label{def:homotopical_presentation:3}
  Let \(C,B\) be categories. Let \(F:C\rightarrow B\) be a
  functor. Let \(\Sigma\) be a set of morphisms in \(C\)
  mapped to an isomorphism in \(B\). For each \(d\in C\), we
  define a subcategory \(A^{\Sigma}_{d}\) of \(C_{d}\) as
  the full subcategory spanned by the morphisms
  \(\phi:e\rightarrow d\) in \(\Sigma\).
\end{definition}

\begin{remark}
  \label{rem:homotopical_presentation:1}
  In Definition~\ref{def:homotopical_presentation:3}, we
  have
  \begin{equation}
    \label{eqn:homotopical_presentation:2}
    \Sigma=S_{A^{\Sigma}}
  \end{equation}
  and, for every family of subcategories \(A_{d}\) of
  \(C_{d}\) for \(d\in C\), \(A_{d}\) is a subcategory of
  \((A^{S_{A}})_{d}\) with the same set of objects for all
  \(d\in C\).
\end{remark}

\begin{lemma}
  \label{lem:homotopical_presentation:2}
  Let \(C\) be a category. Let \(\Sigma\) be a set of
  morphisms in \(C\). Let \(F:C\rightarrow C[\Sigma^{-1}]\)
  be the localization. If \(\Sigma\) admits the calculus of
  right fractions then \(A^{\Sigma}_{d}\) is
  \(F\)-admissible for each \(d\in C\).
\end{lemma}
\begin{proof}
  The localization \(F\) have \(A^{\Sigma}_{d}\)-lifting
  property by
  Lemma~\ref{rem:presheaves_and_Q-divisors:47}(2) and
  \eqref{eqn:homotopical_presentation:2}.
  \(A^{\Sigma}_{d}\) is cofiltered by the dual of
  Proposition~7.1.10 in~\cite{kashiwara-shapira-06}. Thus
  \(A^{\Sigma}_{d}\) is \(F\)-admissible for each
  \(d\in C\).
\end{proof}

\begin{remark}
  \label{rem:presheaves_and_Q-divisors:45}
  Because of the results in this section, the
  \(A_{d}\)-lifting property may be called a right
  \(A_{d}\)-lifting property.
\end{remark}

\section{Homotopical presentations of small categories}
\label{sec:proof-theorem-ref}

In this section, we prove
Theorem~\ref{thm:presheaves_and_Q-divisors:1} and its
corollaries. The proof of
Corollary~\ref{cor:presheaves_and_Q-divisors:1} relies on
Corollary~\ref{cor:presheaves_and_Q-divisors:2}. So the
proof will be completed after we prove the corollary. The
proof for \prfr 1 2 in
Corollary~\ref{cor:presheaves_and_Q-divisors:1} follows from
Theorem~\ref{thm:presheaves_and_Q-divisors:1}.  The proof of
\prfl 1 2 will follow from
Lemma~\ref{lem:presheaves_and_Q-divisors:13},
Lemma~\ref{lem:presheaves_and_Q-divisors:34} and
Corollary~\ref{cor:presheaves_and_Q-divisors:2}.

Corollary~\ref{cor:presheaves_and_Q-divisors:1} is about a special case
of the small presentations of model categories treated in
\cite{dugger-01}, \cite{dugger-01-b}: Given a model category \(M\),
find a small category \(C\) and a set of morphisms \(S\) in \(\Un C\)
such that \(\BL_{S}(\Un C)\) is Quillen equivalent to \(M\).

In our setting \(M\) is \(\Un B\). We are given a functor
\(F:C\rightarrow B\), hence \(F_{*}:\Un C\rightarrow \Un B\).  It is
not always possible to get such a Quillen equivalence from \(F_{*}\)
by localizing at the set of all morphisms mapped to isomorphisms.
For example \(F:\partial\Delta[1]\rightarrow \Delta[0]\) and
\(F:\partial\Delta[1]\rightarrow \Delta[1]\) do not produce such
Quillen equivalences. \(F\) should lift morphisms from \(B\) to \(C\)
so that objects in fibers can be connected when their images are
connected.
% Without it, the total derived functor would not be even essentially
% surjective on objects.
The lifting property
in~Definition~\ref{def:presheaves_and_Q-divisors:2} is
required for this reason.

Theorem~\ref{thm:presheaves_and_Q-divisors:1} follows from
Lemma~\ref{lem:presheaves_and_Q-divisors:5} and
Lemma~\ref{lem:presheaves_and_Q-divisors:7} below.

\begin{lemma}
  \label{lem:presheaves_and_Q-divisors:3}
  Let \(B,C\) be small categories. Let \(F:C\rightarrow B\) be a
  functor. We assume that \(C\) has finite limits and \(F\) preserves
  finite limits.  We assume that the following hold.
  \begin{enumerate}
  \item For each \(d\in C\), there is an \(F\)-admissible
     subcategory
    \(A_{d}\) of \(C_{d}\).
  \item \(F\) is essentially surjective on objects. 
  \end{enumerate}
  Then the counit
  \begin{equation}
    \label{eqn:presheaves_and_Q-divisors:20}
    \varepsilon:F_{*}F^{*}\xrightarrow[\cong]{} I_{\sPSh(B)}
  \end{equation}
  of the adjunction \eqref{eqn:presheaves_and_Q-divisors:168} is an
  isomorphism.
\end{lemma}
\begin{proof}
  We recall from \cite{artin-62} that for each \(X\in \sPSh(C)\) and
  \(b\in B\),
  \begin{equation}
    \label{eqn:presheaves_and_Q-divisors:181}
    F_{*}X(b)=\colim_{\substack{\gangle{c,b\rightarrow Fc}\\\in (b\downarrow F)^{\op}}} X(c).
  \end{equation}
  Let \(Y\) be an object of \(\sPSh(B)\).  Let \(b\) be an object of
  \(B\).  Because \(F\) is essentially surjective on objects, we may
  assume that there is an object \(d\in C\) such that \(b=Fd\).  By
  \eqref{eqn:presheaves_and_Q-divisors:181}, we have
  \begin{equation}
    \label{eqn:presheaves_and_Q-divisors:182}
    F_{*}F^{*}Y(Fd)=\colim_{\substack{\gangle{e,Fd\rightarrow Fe}\\\in (B_{d})^{\op}}} Y(Fe).
  \end{equation}
  Thus, we have the following commutative diagram in \(\sSet\).
  \begin{equation}
    \label{eqn:presheaves_and_Q-divisors:23}
    \begin{tikzcd}
      F_{*}F^{*}Y(Fd) \ar[r,"\varepsilon_{Y,Fd}"]  & Y(Fd) \\
      \colim\limits_{\substack{\gangle{e,e\rightarrow d}\\\in (A_{d})^{\op}}} Y(Fe) \ar[u,"\cong"'] &
      \colim\limits_{\substack{\gangle{e,e\rightarrow d}\\\in (A_{d})^{\op}}}
      Y(Fd) \ar[l,"\cong"] \ar[u,"\cong"']
    \end{tikzcd}
  \end{equation}
  The left vertical map is an isomorphism by
  \eqref{eqn:presheaves_and_Q-divisors:182} and
  Lemma~\ref{lem:presheaves_and_Q-divisors:2}.  The right vertical map
  is an isomorphism because \(A_{d}\) is cofiltered, hence
  connected. The bottom horizontal map is an isomorphism by the
  definition of \(A_{d}\). Thus the counit \(\varepsilon\) is an
  isomorphism.
\end{proof}

\begin{lemma}
  \label{lem:presheaves_and_Q-divisors:4}
  Let \(B,C\) be small categories.  Let \(F:C\rightarrow B\) be a
  functor.  We assume that \(C\) has finite limits and \(F\) preserves
  finite limits.  Let \(X\) be an object of \(\sPSh(C)\). Let
  \begin{equation}
    \label{eqn:presheaves_and_Q-divisors:24}
    \eta_{X}:X\xrightarrow[\simeq]{} F^{*}F_{*}X
  \end{equation}
  be the component at \(X\) of the unit \(\eta\) of the adjunction
  \eqref{eqn:presheaves_and_Q-divisors:168}.  We assume that the
  following hold. 
  \begin{enumerate}
  \item For each \(d\in C\), there is an \(F\)-admissible subcategory
    \(A_{d}\) of \(C_{d}\).
  \item \(X(\phi)\) is a weak equivalence in \(\sSet\) for each
    \(\phi\in S_{A}\).
  \end{enumerate}
  Then for each object \(d\) of \(C\), \(\eta_{X,d}\) is a
  weak equivalence in \(\sSet\).
\end{lemma}
\begin{proof}
  Let \(d\) be an object of \(C\).
  We have the following commutative diagram in \(\sSet\).
  \begin{equation}
    \label{eqn:presheaves_and_Q-divisors:26}
    \begin{tikzcd}
      X(d) \ar[r,"\eta_{X,d}"] \ar[d,"="] & F^{*}F_{*}X(d) \ar[r,"="'] & F_{*}X(Fd) \\
      X(d) & \colim\limits_{\substack{\gangle{e,e\rightarrow d}\\\in
          (A_{d})^{\op}}} X(d) \ar[l,"\cong"] \ar[r,"\simeq"'] &
      \colim\limits_{\substack{\gangle{e,e\rightarrow d}\\\in
          (A_{d})^{\op}}} X(e) \ar[u,"\cong"']
    \end{tikzcd}
  \end{equation}
  The right vertical map is an isomorphism by
  \eqref{eqn:presheaves_and_Q-divisors:181} and
  Lemma~\ref{lem:presheaves_and_Q-divisors:2}.  The left horizontal
  map in the bottom is an isomorphism because \(A_{d}\) is connected.
  The right horizontal map in the bottom is a weak equivalence because
  \((A_{d})^{\op}\) is filtered and \(X(\phi)\) is a weak equivalence
  in \(\sSet\) for each \(\phi\in S_{A}\) (cf.~Theorem~14.10
  in~\cite{hirschhorn-14}). Thus \(\eta_{X,d}\) is a weak equivalence
  in \(\sSet\).
\end{proof}

\begin{lemma}
  \label{lem:presheaves_and_Q-divisors:5}
  Let \(B,C\) be small categories. Let \(F:C\rightarrow B\) be a
  functor. We assume that \(C\) has finite limits and \(F\) preserves
  finite limits. We assume that the following hold.
  \begin{enumerate}
  \item For each \(d\in C\), there is an \(F\)-admissible subcategory
    \(A_{d}\) of \(C_{d}\).
  \end{enumerate}
  Let
  \begin{equation}
    \label{eqn:presheaves_and_Q-divisors:85}
    F_{*}:\BL_{S_{A}}(\Un C)\rightleftarrows \Un B:F^{*}
  \end{equation}
  be the Quillen adjunction obtained from the universal adjunction of
  \(F\).  Let \(f:X\rightarrow Y\) be a morphism between cofibrant
  objects of \(\BL_{S_{A}}(\Un C)\). If \(F_{*}f\) is a weak
  equivalence in \(\Un B\) then \(f\) is a weak equivalence in
  \(\BL_{S_{A}}(\Un C)\).
\end{lemma}
\begin{proof}
  Consider the commutative diagram
  \begin{equation}
    \label{eqn:presheaves_and_Q-divisors:27}
    \begin{tikzcd}
      X \ar[r, "f"] \ar[d,"i_{X}"',"\simeq"] & Y \ar[d,"i_{Y}","\simeq"']\\
      RX \ar[r,"Rf"]& RY
    \end{tikzcd}
  \end{equation}
  where \(i_{X}\), \(i_{Y}\) are fibrant replacements of
  \(X\), \(Y\) in \(\BL_{S_{A}}(\Un C)\) respectively.
  Morphisms \(F_{*}i_{X}\) and \(F_{*}i_{Y}\) are weak
  equivalences in \(\Un B\).  Thus \(F_{*}Rf\) is also a
  weak equivalence in \(\Un B\).  Consider the following
  commutative diagram.
  \begin{equation}
    \label{eqn:presheaves_and_Q-divisors:28}
    \begin{tikzcd}[column sep=large]
      RX \ar[r,"Rf"] \ar[d,"\eta_{RX}"',"\simeq"] & RY \ar[d,"\eta_{RY}","\simeq"']\\
      F^{*}F_{*}RX \ar[r,"F^{*}F_{*}Rf","\simeq"'] & F^{*}F_{*}RY
    \end{tikzcd}
  \end{equation}
  The vertical morphisms are weak equivalences in \(\Un C\)
  by Lemma~\ref{lem:presheaves_and_Q-divisors:6} and
  Lemma~\ref{lem:presheaves_and_Q-divisors:4} because \(RX\)
  and \(RY\) are fibrant objects of \(\BL_{S_{A}}(\Un C)\).
  \(F^{*}F_{*}Rf\) is a weak equivalence in \(\Un C\)
  because \(F_{*}Rf\) is a weak equivalence in \(\Un
  B\). Thus \(f\) is a weak equivalence in
  \(\BL_{S_{A}}(\Un C)\).
\end{proof}

\begin{lemma}
  \label{lem:presheaves_and_Q-divisors:7}
  Let \(B,C\) be small categories. Let \(F:C\rightarrow B\) be a
  functor. We assume that \(C\) has finite limits and \(F\) preserves
  finite limits.  We assume that the following hold.
  \begin{enumerate}
  \item\label{item:16} For each \(d\in C\), there is an \(F\)-admissible
     subcategory
    \(A_{d}\) of \(C_{d}\).
  \item\label{item:17} \(F\) is essentially surjective on objects. 
  \end{enumerate}
  Let
  \begin{equation}
    \label{eqn:presheaves_and_Q-divisors:19}
    F_{*}:\BL_{S_{A}}(\Un C)\rightleftarrows \Un B:F^{*}
  \end{equation}
  be the Quillen adjunction obtained from the universal adjunction of
  \(F\).  Then, for each object \(Y\) of \(\Un B\), the composite
  \begin{equation}
    \label{eqn:presheaves_and_Q-divisors:29}
    F_{*}QF^{*}Y\xrightarrow{F_{*}(p_{F^{*}Y})}
    F_{*}F^{*}Y\xrightarrow{\varepsilon_{Y}} Y
  \end{equation}
  is a weak equivalence in \(\Un B\) where \(p_{F^{*}Y}\) is
  a cofibrant replacement of \(F^{*}Y\) in
  \(\BL_{S_{A}}(\Un C)\) and \(\varepsilon_{Y}\) is the
  component at \(Y\) of the counit \(\varepsilon\) of
  \eqref{eqn:presheaves_and_Q-divisors:168}.
\end{lemma}
\begin{proof}
  For each \(b\in B\), \(F_{*}(p_{F^{*}Y})_{b}\) is computed over
  \((b\downarrow F)^{\op}\).  Because \((b\downarrow F)^{\op}\) is
  filtered for each \(b\in B\) by
  Lemma~\ref{lem:presheaves_and_Q-divisors:30} and \(p_{F^{*}Y}\) is a
  weak equivalence in \(\Un C\), \(F_{*}(p_{F^{*}Y})\) is a weak
  equivalence in \(\Un B\). The component \(\varepsilon_{Y}\) is an
  isomorphism by Lemma~\ref{lem:presheaves_and_Q-divisors:3}. Thus the
  composite is a weak equivalence in \(\Un B\).
\end{proof}

\begin{proof}[Proof of Theorem~\ref{thm:presheaves_and_Q-divisors:1}]
  Because the elements of \(S_{A}\) are mapped to
  isomorphisms in \(B\), the universal adjunction of \(F\)
  extends to the Quillen adjunction
  \eqref{eqn:presheaves_and_Q-divisors:135}.  So, it follows
  from Corollary~1.3.16 in~\cite{hovey-99},
  Lemma~\ref{lem:presheaves_and_Q-divisors:5} and
  Lemma~\ref{lem:presheaves_and_Q-divisors:7}.
\end{proof}

The following is a partial converse of
Theorem~\ref{thm:presheaves_and_Q-divisors:1}. 

\begin{lemma}
  \label{lem:presheaves_and_Q-divisors:13}
  Let \(B,C\) be small categories. Let \(F:C\rightarrow B\)
  be a functor. We assume that \(C\) has finite limits and
  \(F\) preserves finite limits. We assume that we are give
  a family of subcategory \(A_{d}\) of \(C_{d}\) for each
  \(d\in C\). Let
  \begin{equation}
    \label{eqn:presheaves_and_Q-divisors:67}
    F_{*}:\BL_{S_{A}}(\Un C)\rightleftarrows \Un B:F^{*}
  \end{equation}
  be the Quillen adjunction obtained from the universal adjunction of
  \(F\).  We also assume that the following hold.
  \begin{enumerate}
  \item\label{item:14} For each \(d\in C\), \(F\) has
    \(A_{d}\)-lifting property.
  \item The adjunction \eqref{eqn:presheaves_and_Q-divisors:67} is a
    Quillen equivalence.
  \end{enumerate}
  Then \(F\) is essentially surjective on objects.
\end{lemma}
\begin{proof}
  Let \(Y\) be an object of \(\Un B\).
  Consider a commutative diagram
  \begin{equation}
    \label{eqn:presheaves_and_Q-divisors:227}
    \begin{tikzcd}
      [column sep=2cm]
      F_{*}QF^{*}Y \ar[r,"{F_{*}(p_{F^{*}Y})}"]
      \ar[d,"F_{*}QF^{*}i_{Y}"] &
      F_{*}F^{*}Y\ar[r,"{\varepsilon_{Y}}"] \ar[d,"F_{*}F^{*}i_{Y}"] &
      Y \ar[d,"i_{Y}"] \\
      F_{*}QF^{*}RY \ar[r,"{F_{*}(p_{F^{*}RY})}"] &
      F_{*}F^{*}RY\ar[r,"{\varepsilon_{RY}}"]& RY
    \end{tikzcd}
  \end{equation}
  where \(p_{F^{*}Y}\) is a cofibrant replacement of \(F^{*}Y\) in
  \(\BL_{S_{A}}(\Un C)\) and \(i_{Y}\) is a fibrant replacement of
  \(Y\) in \(\Un B\). The morphism \(QF^{*}i_{Y}\) is a weak
  equivalence in \(\BL_{S_{A}}(\Un C)\) between cofibrant objects.
  So, \(F_{*}QF^{*}i_{Y}\) is a weak equivalence in \(\Un B\). Then,
  because \eqref{eqn:presheaves_and_Q-divisors:67} is a Quillen
  equivalence, the top composite is a weak equivalence in \(\Un B\).
  Since \((b\downarrow F)^{\op}\) is filtered for each \(b\in B\) by
  Lemma~\ref{lem:presheaves_and_Q-divisors:30}, the top left morphism
  is a weak equivalence in \(\Un B\).  Thus, \(\varepsilon_{Y}\) is
  also a weak equivalence. Then for each \(b\in B\),
  \begin{equation}
    \label{eqn:presheaves_and_Q-divisors:228}
    \varepsilon_{Y,b}:
    \colim\limits_{\substack{\gangle{c,b\rightarrow Fc}\\\in
        (b\downarrow F)^{\op}}}
    Y(Fc)\rightarrow Y(b)    
  \end{equation}
  is a weak equivalence in \(\sSet\).  

  Let \(b\) be an object of \(B\). We apply the weak equivalence
  \eqref{eqn:presheaves_and_Q-divisors:228} to \(r_{B}(b)\). Let
  \(Y=r_{B}(b)\).  Because \(\hom_{B}(b,b)=Y_{0}(b)=\pi_{0}(Y(b))\),
  \eqref{eqn:presheaves_and_Q-divisors:228} induces a bijection
  \begin{equation}
    \label{eqn:presheaves_and_Q-divisors:229}
    \colim\limits_{\substack{\gangle{c,b\rightarrow Fc}\\\in
        (b\downarrow F)^{\op}}}
    \hom_{B}(Fc,b)
    \xrightarrow\cong \hom_{B}(b,b).
  \end{equation}
  Thus there is an object \(d\in C\) and morphisms
  \(\alpha:b\rightarrow Fd\) and \(\beta:Fd\rightarrow b\) such that
  \begin{equation}
    \label{eqn:presheaves_and_Q-divisors:230}
    \beta\cdot\alpha=\id_{b}.
  \end{equation}
  Then \(\alpha\cdot\beta:Fd\rightarrow Fd\) is a split
  idempotent in \(B\) by
  \eqref{eqn:presheaves_and_Q-divisors:230}. Thus
  \(\alpha:b\rightarrow Fd\) is the equalizer of
  \(\alpha\cdot \beta\) and \(1_{Fd}\) by Proposition~1
  in~\cite{borceux-dejean-86}.
  \begin{equation}
    \label{eqn:presheaves_and_Q-divisors:233}
    \alpha= \Eq\left(
    \begin{tikzcd}
      Fd \ar[r,shift left, "{\alpha\cdot\beta}"] \ar[r, shift right,
      "{1_{Fd}}"']
      & Fd
    \end{tikzcd}
    \right)
  \end{equation}

  By the assumption \eqref{item:14}, there are morphisms
  \(\phi,\psi:e\rightarrow d\) in \(C\) such that \(F\phi\) is an
  isomorphism and
  \begin{align}
    \label{eqn:presheaves_and_Q-divisors:231}
    F\psi&=\alpha\cdot\beta\cdot F\phi\\
    F\phi&=1_{Fd}\cdot F\phi
  \end{align}
  hold. Let \(\varepsilon:e'\rightarrow e\) be the equalizer
  of \(\phi,\psi\) in \(C\).  Since \(F\) preserves finite
  limits, \(F\varepsilon:Fe'\rightarrow Fe\) is the
  equalizer of \(F\phi\) and \(F\psi\) in \(B\). Then, since
  \(F\phi\) is an isomorphism,
  \(F\phi\cdot F\varepsilon:Fe'\rightarrow Fd\) is the
  equalizer of \(\alpha\cdot\beta\) and \(1_{Fd}\).
  \begin{equation}
    \label{eqn:presheaves_and_Q-divisors:232}
    F\phi\cdot F\varepsilon= \Eq\left(
    \begin{tikzcd}
      Fd \ar[r,shift left, "{\alpha\cdot\beta}"] \ar[r, shift right,
      "{1_{Fd}}"']
      & Fd
    \end{tikzcd}
    \right)
  \end{equation}
  Thus \(Fe'\cong b\) by \eqref{eqn:presheaves_and_Q-divisors:233} and
  \eqref{eqn:presheaves_and_Q-divisors:232}, and \(F\) is essentially
  surjective on objects.
\end{proof}

We prepare a lemma for the proof of
Corollary~\ref{cor:presheaves_and_Q-divisors:6}.

\begin{lemma}
  \label{lem:presheaves_and_Q-divisors:34}
  Let \(B,C\) be small categories. Let \(F:C\rightarrow B\) a functor.
  If \(F\) is a universal equivalence then \(F\) is full and faithful.
\end{lemma}
\begin{proof}
  Let \(e\) be an object of \(C\). Let \(X=r_{C}(e)\).  By our
  assumption the composite
  \begin{equation}
    \label{eqn:presheaves_and_Q-divisors:242}
    X\xrightarrow[\simeq]{\eta_{X}}
    F^{*}F_{*}X\xrightarrow[\simeq]{F^{*}(i_{F_{*}X})} F^{*}RF_{*}X
  \end{equation}
  is a weak equivalence in \(\Un C\) where \(\eta\) is the unit of
  \eqref{eqn:presheaves_and_Q-divisors:168}. Because \(F^{*}\)
  preserves weak equivalences, \(F^{*}(i_{F_{*}X})\) is a weak
  equivalence in \(\Un C\). So, the morphism \(\eta_{X}\) is a weak
  equivalence \(\Un C\).  Then for each \(d\in C\), we have an
  isomorphism
  \begin{equation}
    \label{eqn:presheaves_and_Q-divisors:243}
    \hom_{C}(d,e)\xrightarrow[\cong]{F} \hom_{B}(Fd,Fe)
  \end{equation}
  by Lemma~\ref{lem:presheaves_and_Q-divisors:20}.
  Thus \(F\) is full and faithful.
\end{proof}

\begin{proof}[Proof of
  Corollary~\ref{cor:presheaves_and_Q-divisors:6}]
  First, we make an observation.  If \(F\) is full and faithful then
  \begin{enumerate}
  \item[(i)] \(S_{C}\) is the set of all isomorphisms in \(C\), hence
    \(\BL_{S_{C}}(\Un C)=\Un C\), and
  \item[(ii)] \(F\) has \(C_{d}\)-lifting property for all
    \(d\in C\) for a trivial reason.
  \end{enumerate}
  If either (1) or (2) holds then \(F\) is full and faithful by
  Lemma~\ref{lem:presheaves_and_Q-divisors:34}. Thus
  Corollary~\ref{cor:presheaves_and_Q-divisors:1} reduces to
  Corollary~\ref{cor:presheaves_and_Q-divisors:6}.
%  First, we make an observation.  If \(F\) is full and faithful then
%  \begin{enumerate}
%  \item[(i)] the set \(S_{C}\) is the set of all isomorphisms in
%    \(C\), hence \(\BL_{S_{C}}(\Un C)=\Un C\), and
%  \item[(ii)] \(F\) satisfies (Lift) for a trivial reason.
%  \end{enumerate}

%  \prfr 1 2 If \(F\) is an equivalence of categories then \(F\) is
%  essentially surjective on objects. So, by (i), (ii) above and
%  Corollary~\ref{cor:presheaves_and_Q-divisors:1}, \(F\) is a universal
%  equivalence.

%  \prfl 1 2 First, \(F\) is full and faithful by
%  Lemma~\ref{lem:presheaves_and_Q-divisors:34}.  Then (i) above and
%  Corollary~\ref{cor:presheaves_and_Q-divisors:1} imply that \(F\) is
%  essentially surjective on objects.
\end{proof}

\begin{proof}[Proof of
  Corollary~\ref{cor:presheaves_and_Q-divisors:2}]
  We apply the part \prfr 1 2 of
  Corollary~\ref{cor:presheaves_and_Q-divisors:1} to the localization
  \(F:C\rightarrow \Sigma^{-1}C\) of \(C\) at \(\Sigma\).  So,
  \(S_{C}\) is the set of all morphisms in \(C\) mapped to an
  isomorphism in \(\Sigma^{-1}C\) by \(F\).  We have
  \(\Sigma\subseteq S_{C}\).  Because \(\Sigma\) admits the calculus
  of right fractions,
  \begin{enumerate}
  \item the localization \(F:C\rightarrow \Sigma^{-1}C\) preserves
    finite limits by the dual of Proposition on Chapter I.3.1
    in~\cite{gabriel-zisman-67}, and
  \item \(F\) has \(C_{d}\)-lifting property for all \(d\in C\).
  \end{enumerate}
  So we have a Quillen equivalence
  \begin{equation}
    \label{eqn:presheaves_and_Q-divisors:248}
    F_{*}:\BL_{S_{C}}(\Un C)\rightleftarrows \Un{\Sigma^{-1}C}:F^{*}
  \end{equation}
  by \prfr 1 2 of Corollary~\ref{cor:presheaves_and_Q-divisors:1}.  As
  the set of weak equivalences in \(\sSet\) satisfies the 2-out-of-6
  property, the \(r_{C}(S_{C})\)-local objects are precisely the
  \(r_{C}(\Sigma)\)-local objects by the dual of the description of
  \(S_{C}\) on Chapter~I.3.5 in~\cite{gabriel-zisman-67}
  (Proposition~7.1.20 in~\cite{kashiwara-shapira-06}). Thus
  \begin{equation}
    \label{eqn:presheaves_and_Q-divisors:249}
    \BL_{\Sigma}(\Un C)=\BL_{S_{C}}(\Un C)
  \end{equation}
  holds. Thus \eqref{eqn:presheaves_and_Q-divisors:106} is a Quillen
  equivalence.
\end{proof}

\begin{remark}
  \label{rem:presheaves_and_Q-divisors:25}
  Instead of Corollary~\ref{cor:presheaves_and_Q-divisors:1},
  we can use Theorem~\ref{thm:presheaves_and_Q-divisors:1}
  and Lemma~\ref{lem:homotopical_presentation:2} to prove
  Corollary~\ref{cor:presheaves_and_Q-divisors:2} directly.
\end{remark}

\begin{proof}[Proof of Corollary~\ref{cor:presheaves_and_Q-divisors:1}]
  \prfr 1 2 follows from
  Theorem~\ref{thm:presheaves_and_Q-divisors:1}. Now, we
  consider \prfl 1 2.  Because of
  Lemma~\ref{lem:presheaves_and_Q-divisors:13}, it is enough
  to show that \(F\) has \(C_{d}\)-lifting property for all
  \(d\in C\).  Because \(C\) has finite limits and \(F\)
  preserves finite limits, \(S_{C}\) admits the calculus of
  right fractions.  Then, the functor
  \(H:C[\left(S_{C}\right)^{-1}]\rightarrow B\) associated
  with \(F\) is a universal equivalence by (2) and
  Corollaries~\ref{cor:presheaves_and_Q-divisors:2}.
  % The category \(S_{C}^{-1}C\) has finite limits and the functor
  % \(H\) preservse finite limits by the dual of Corollary~1 on
  % Chapter~I.3.2 in \cite{gabriel-zisman-67} and its proof.
  Then \(H\) is full by Lemma~\ref{lem:presheaves_and_Q-divisors:34}.
  % by Corollary~\ref{cor:presheaves_and_Q-divisors:6}.
  Thus \(F\) has \(C_{d}\)-lifting property for all
  \(d\in C\) by Lemma~\ref{lem:homotopical_presentation:1}.
\end{proof}

\begin{remark}
  \label{rem:presheaves_and_Q-divisors:30}
  Let \(B,C\) be small categories. Let \(F:C\rightarrow B\) be a
  functor. We assume that \(C\) has finite limits and \(F\) preserves
  finite limits. Then by an argument similar to the proof of \prfl 1 2
  in Corollary~\ref{cor:presheaves_and_Q-divisors:1}, we can show that
  the following are equivalent.
  \begin{enumerate}
  \item \(F\) generates a homotopical presentation of \(B\)
    with a generator \(C\) and a relation \(S_{C}\).
  \item There is a subset \(\Sigma\) of \(S_{C}\) admitting
    the calculus of right fraction such that \(F\) generates
    a homotopical pretension of \(B\) with a generator \(C\)
    and a relation \(\Sigma\).
%  \item For each \(d\in C\), there is an \(F\)-admissible subcategory
%    \(A_{d}\) of \(C_{d}\) such that 
%    \begin{equation}
%      \label{eqn:presheaves_and_Q-divisors:285}
%      F_{*}:\BL_{S_{A}}(\Un C)\rightleftarrows \Un B:F^{*}
%    \end{equation}
%    is a Quillen equivalence.    
  \end{enumerate}
  The point of \prfl 1 2 is that we can localizer further up
  to \(S_{C}\) keeping the Quillen equivalence.
\end{remark}

\section{Examples}
\label{sec:examples}

This section contains some examples of homotopical presentations.

\begin{example}
  [Chapter~I.2.5 in~\cite{gabriel-zisman-67}]  
  \label{exa:presheaves_and_Q-divisors:5}
  Let \(A\) be a small abelian category. Let \(B\) be a thick
  subcategory of \(A\). Let \(\Sigma=\Sigma(B)\) be the set of all
  morphisms \(f\) in \(A\) such that \(\ker(f)\) and \(\coker(f)\)
  belong to \(B\). Then \(\Sigma\) satisfies the calculus of left and
  right fractions. The category \(\Sigma^{-1}A\) is isomorphic to the
  quotient category \(A/B\) on Chapter~1.11
  in~\cite{grothendieck}. Thus we have a Quillen equivalence
  \begin{equation}
    \label{eqn:presheaves_and_Q-divisors:250}
    \BL_{\Sigma}(\Un A)\rightleftarrows \Un{A/B}
  \end{equation}
  by Corollary~\ref{cor:presheaves_and_Q-divisors:2}.
\end{example}

% \begin{example}
%   \label{exa:presheaves_and_Q-divisors:6}
%   Let \(A\) be an abelian category. Let \(\Sigma\) be the set of
%   essential extensions in \(A\). Then \(\Sigma\) satisfies the
%   calculus of right fractions.
% \end{example}

\begin{remark}
  \label{rem:presheaves_and_Q-divisors:21}
  One cannot apply Corollary~\ref{cor:presheaves_and_Q-divisors:2} to
  the categories of fibrant objects (\cite{brown-73}) in general. Let
  \((M;W,C)\) be a category of fibrant objects.  Brown defined a
  quotient category \(\pi M\) as an approximation to the homotopy
  category of \(M\) (p.423 in~\cite{brown-73}).  While the image of
  \(W\) in \(\pi M\) satisfies the calculus of right fractions by
  Proposition~2 in Part I.2 in~\cite{brown-73}, \(\pi M\) does not
  have finite limits in general.  If it does, then the homotopy
  category \(\Ho M\) would have finite limits.  For example, every
  model category has an associated category of fibrant objects whose
  homotopy category is equivalent to that of the model category.  But
  the homotopy category of a model category does not have finite
  limits in general.
\end{remark}

We will use the following proposition to provide examples of
functors generating homotopical presentations. It follows
from Corollary~\ref{cor:presheaves_and_Q-divisors:6} and
Corollary~\ref{cor:presheaves_and_Q-divisors:2}.

\begin{proposition}
  \label{pro:presheaves_and_Q-divisors:1}
  Let \(F:C\rightleftarrows B:G\) be an adjunction between
  small categories.  We assume that the following hold.
  \begin{enumerate}
  \item \(C\) has finite limits.
  \item \(S_{C}\) admits the calculus of right fractions.
  \item The right adjoint \(G\) is fully faithful.
  \end{enumerate}
  Let \(P:C\rightarrow(S_{C})^{-1}C\) be the localization of \(C\) at
  \(S_{C}\) and \(H:(S_{C})^{-1}C\rightarrow B\) be the functor such
  that \(F=H\cdot P\).  Then
  \begin{equation}
    \label{eqn:presheaves_and_Q-divisors:269}
    P_{*}:\BL_{S_{C}}(\Un{C})
    \rightleftarrows
    \Un{(S_{C})^{-1}C}:P^{*}
  \end{equation}
  and
  \begin{equation}
    \label{eqn:presheaves_and_Q-divisors:251}
    H_{*}:\Un{(S_{C})^{-1}C}
    \rightleftarrows
    \Un{B}:H^{*}    
  \end{equation}
  are Quillen equivalences.  Thus, \(F\) generates a homotopical
  presentation of \(B\).
\end{proposition}
\begin{proof}
  By (1), (2) and Corollary~\ref{cor:presheaves_and_Q-divisors:2}, we
  have the Quillen equivalence
  \eqref{eqn:presheaves_and_Q-divisors:269}.  The category
  \((S_{C})^{-1}C\) has finite limits by the dual of Corollary~1 on
  Chapter~I.3.2 in~\cite{gabriel-zisman-67} (Proposition~7.1.22
  in~\cite{kashiwara-shapira-06}). By (3) and the dual of Proposition
  on Chapter I.1.3 in~\cite{gabriel-zisman-67}, \(H\) is an
  equivalence of categories.  So, we have the Quillen equivalence
  \eqref{eqn:presheaves_and_Q-divisors:251} by
  Corollary~\ref{cor:presheaves_and_Q-divisors:6}. Thus
  \begin{equation}
    \label{eqn:presheaves_and_Q-divisors:271}
    F_{*}:\BL_{S_{C}}(\Un C)\rightleftarrows \Un B:F^{*}
  \end{equation}  
  is a Quillen equivalence.
\end{proof}

\begin{remark}
  \label{rem:presheaves_and_Q-divisors:31}
  In Proposition~\ref{pro:presheaves_and_Q-divisors:1}, (1) implies
  (2) if \(F\) preserves finite limits.
\end{remark}

In the following examples, we will implicitly assume that we
fixed a universe in each example relative to which the
categories in the examples are small.

In Example~\ref{exa:presheaves_and_Q-divisors:4}, we show that the
associated sheaf functor generates a homotopical presentation of the
category of sheaves with the set of local isomorphisms as the set of
relations.
  
\begin{example}
  \label{exa:presheaves_and_Q-divisors:4}
  Let \(T\) be a site. Let \(\PSh(T)\) be the category of presheaves
  on (the underlying category of) \(T\) and \(\Sh(T)\) be the category
  of sheaves on \(T\).  We have an adjunction
  \begin{equation}
    \label{eqn:presheaves_and_Q-divisors:272}
    \begin{tikzcd}
      \PSh(T) \ar[rr,shift left=2,"a"] &&
      \Sh(T) \ar[ll,shift left=2,"\bot"',"i"]
    \end{tikzcd}    
  \end{equation}
  where \(a\) is the associated sheaf functor and \(i\) is
  the embedding.  Let \(\LI\) be the set of local
  isomorphisms in \(\PSh(T)\).  A morphism
  \(f:X\rightarrow Y\) in \(\PSh(T)\) is a local
  isomorphisms iff \(a(f):a(X)\rightarrow a(Y)\) is an
  isomorphisms in \(\Sh(T)\) (Exercise~16.7
  in~\cite{kashiwara-shapira-06}). Since the functor \(a\)
  preserves finite limits, \(\LI\) admits the calculus of
  right fractions (Exercise~7.5
  in~\cite{kashiwara-shapira-06}). So, if we fix a universe
  relative to which \(\PSh(T)\) and \(\Sh(T)\) are small,
  the associated sheaf functor
  \(a:\PSh(T)\rightarrow \Sh(T)\) generates a homotopical
  presentation of \(\Sh(T)\) with \(\LI\) as the set of
  relations
  \begin{equation}
    \label{eqn:presheaves_and_Q-divisors:268}
    a_{*}:\BL_{\LI}(\Un{\PSh(T)})
    \rightleftarrows
    \Un{\Sh(T)}:a^{*}
  \end{equation}  
  by Proposition~\ref{pro:presheaves_and_Q-divisors:1}.
\end{example}

In Example~\ref{exa:presheaves_and_Q-divisors:2}, we show that, in
general, \(\Un -\) does not preserve Quillen equivalences.

\begin{example}
  \label{exa:presheaves_and_Q-divisors:2}
  Consider the Quillen equivalence
  \(|\cdot |:\sSet\rightleftarrows\Top:Sing\) associated
  with the geometric realization functor and the singular
  simplicial set functor. Let \(\Sigma\) be the set of all
  morphisms in \(\sSet\) mapped to an homeomorphisms in
  \(\Top\) by \(|\cdot|\).  Since \(|\cdot |\) preserves
  finite limits \(\Sigma\) admits the calculus of right
  fractions. \(Sing\) is fully faithful. So, we have a
  Quillen equivalence
  \begin{equation}
    \label{eqn:presheaves_and_Q-divisors:267}
    |\cdot|_{*}:
    \BL_{\Sigma}(\Un{\sSet})\rightleftarrows\Un{\Top}
    :|\cdot|^{*}
  \end{equation}
  by Proposition~\ref{pro:presheaves_and_Q-divisors:1}.  Now, suppose
  that we also have a Quillen equivalence
  \begin{equation}
    \label{eqn:presheaves_and_Q-divisors:103}
    |\cdot|_{*}:
    \Un{\sSet}\rightleftarrows\Un{\Top}
    :|\cdot|^{*}.
  \end{equation}
  Then
  \begin{equation}
    \label{eqn:presheaves_and_Q-divisors:275}
    \id:
    \Un{\sSet}\rightleftarrows \BL_{\Sigma}(\Un{\sSet})
    :\id
  \end{equation}
  is a Quillen equivalence. Hence,
  \begin{equation}
    \label{eqn:presheaves_and_Q-divisors:276}
    P_{*}:
    \Un{\sSet}\rightleftarrows \Un{\Sigma^{-1}\sSet}
    :P^{*}
  \end{equation}
  is a Quillen equivalence by
  Corollary~\ref{cor:presheaves_and_Q-divisors:2} where
  \(P:\sSet\rightarrow \Sigma^{-1}\sSet\) is the localization.  Thus
  \(P\) is full and faithful by
  Lemma~\ref{lem:presheaves_and_Q-divisors:34}.
  %an equivalence of categories by
  %Corollary~\ref{cor:presheaves_and_Q-divisors:6} because \(P\)
  %preserves finite limits by Proposition I.3.1
  %in~\cite{gabriel-zisman-67}.
  But \(P\) is not faithful.  For
  example,
  \begin{equation}
    \label{eqn:presheaves_and_Q-divisors:277}
    \sSet(\Delta^{0},\Delta^{1})
    \xrightarrow{P}
    \Sigma^{-1}\sSet(\Delta^{0},\Delta^{1})
  \end{equation}
  is not injective.  Thus \eqref{eqn:presheaves_and_Q-divisors:103} is
  not a Quillen equivalence and \(\Un -\) does not preserve Quillen
  equivalence in general.
\end{example}
  
\begin{remark}
  \label{rem:presheaves_and_Q-divisors:23}
  There is the dual notion for the universal model category (Section~4
  in~\cite{dugger-01}). Let \(C\) be a small category. Let \(\Vn C\)
  be the model category whose underlying category is
  \(\left(\sSet^{C}\right)^{\op}=\left(\Un{C^{\op}}\right)^{\op}\) and
  whose model structure is the opposite of \(\Un{C^{\op}}\). It is
  called a co-universal model category of \(C\). For example, the dual
  of Corollary~\ref{cor:presheaves_and_Q-divisors:2} is that, for each
  category \(C\) with finite colimits and a subset \(\Sigma\) of
  \(\mor C\) admitting the calculus of left fractions, there is a
  Quillen equivalence
  \begin{equation}
    \label{eqn:presheaves_and_Q-divisors:257}
    \left((F^{\op})^{*}\right)^{\op}:
    \Vn{\Sigma^{-1}C}\rightleftarrows \BR_{\Sigma}(\Vn C):
    \left((F^{\op})_{*}\right)^{\op}
  \end{equation}
  where \(\BR_{\Sigma}\) is the right Bousfield localization at
  \(\Sigma\). 
\end{remark}

\begin{remark}
  \label{rem:presheaves_and_Q-divisors:18}
  For simplicity, given a functor \(F:C\rightarrow B\), we denote
  \(\left((F^{\op})_{*}\right)^{\op}\) and
  \(\left((F^{\op})^{*}\right)^{\op}\) by \(F_{*}\) and \(F^{*}\)
  respectively.
\end{remark}

\begin{proposition}
  \label{pro:presheaves_and_Q-divisors:2}
  Let \(F:C\rightleftarrows B:G\) be an adjunction between
  small categories.  We assume that the following hold.
  \begin{enumerate}
  \item \(C\) has finite colimits.
  \item The right adjoint \(G\) is fully faithful.
  \end{enumerate}
  Then there is a Quillen equivalence
  \begin{equation}
    \label{eqn:presheaves_and_Q-divisors:274}
    F^{*}:\Vn{B}\rightleftarrows \BR_{S_{C}}\left({\Vn{C}}\right):F_{*}.
  \end{equation}
  where \(\BR_{S_{C}}\left({\Vn{C}}\right)\) is the right Bousfield
  localization at \(S_{C}\).  
\end{proposition}
\begin{proof}
  The proof is similar to that of
  Proposition~\ref{pro:presheaves_and_Q-divisors:1}. Let
  \(P:C\rightarrow(S_{C})^{-1}C\) be the localization of \(C\) at
  \(S_{C}\) and \(H:(S_{C})^{-1}C\rightarrow B\) be the functor such
  that \(F=H\cdot P\). The set \(S_{C}\) admits the calculus of left
  fractions by (1) and Exercise~7.5 in~\cite{kashiwara-shapira-06}.
  Then, by (1) and the dual of
  Corollary~\ref{cor:presheaves_and_Q-divisors:2}, we have the
  following Quillen equivalence.
  \begin{equation}
    \label{eqn:presheaves_and_Q-divisors:273}
    P^{*}:\Vn{(S_{C})^{-1}C}
    \rightleftarrows
    \BR_{S_{C}}(\Vn{C})
    :P_{*}
  \end{equation}
  The category \((S_{C})^{-1}C\) has finite colimits by Corollary~1 on
  Chapter~I.3.2 in~\cite{gabriel-zisman-67} (Proposition~7.1.22
  in~\cite{kashiwara-shapira-06}). By (2) and Proposition on Chapter
  I.1.3 in~\cite{gabriel-zisman-67}, \(H\) is an equivalence of
  categories.  So, we have the Quillen equivalence
  \begin{equation}
    \label{eqn:presheaves_and_Q-divisors:270}
    H^{*}:    \Vn{B}    \rightleftarrows
    \Vn{(S_{C})^{-1}C}
    :H_{*}    
  \end{equation}
  by the dual of Corollary~\ref{cor:presheaves_and_Q-divisors:6}. Thus
  \eqref{eqn:presheaves_and_Q-divisors:274} is a Quillen equivalence.
\end{proof}

\begin{example}
  \label{exa:presheaves_and_Q-divisors:8}
  The following are three examples
  Proposition~\ref{pro:presheaves_and_Q-divisors:2}
  provides.  In
  Example~\ref{exa:presheaves_and_Q-divisors:4}, we have the
  Quillen equivalence
  \begin{equation}
    \label{eqn:presheaves_and_Q-divisors:202}
    a^{*}:\Vn{\Sh(T)}\rightleftarrows \BR_{\LI}{\Vn{\PSh(T)}}:a_{*}.
  \end{equation}
  In Example~\ref{exa:presheaves_and_Q-divisors:2}, we have the 
  Quillen equivalence
  \begin{equation}
    \label{eqn:presheaves_and_Q-divisors:266}
    |\cdot|^{*}:\Vn{\Top}\rightleftarrows
    \BR_{\Sigma}(\Vn{\sSet}):|\cdot|_{*}.    
  \end{equation}
  Consider the adjunction \(h:\sSet\rightleftarrows\Cat:N\) associated
  with the nerve functor \(N\). Since the nerve functor \(N\) is fully
  faithful, we have a Quillen equivalence
  \begin{equation}
    \label{eqn:presheaves_and_Q-divisors:265}
    h^{*}:\Vn{\Cat}\rightleftarrows \BR_{\Sigma}(\Vn{\sSet}):h_{*}    
  \end{equation}
  where \(\Sigma\) is the set of all morphisms in \(\sSet\) mapped to
  an isomorphisms in \(\Cat\) by \(h\).
\end{example}

\begin{remark}
  \label{rem:presheaves_and_Q-divisors:27}
  Theorem~\ref{thm:presheaves_and_Q-divisors:1} and
  Corollary~\ref{cor:presheaves_and_Q-divisors:1} have
  applications in algebraic geometry. The ampleness is one
  of the central notions for positivity of varieties.  One
  can use Corollary~\ref{cor:presheaves_and_Q-divisors:1} to
  characterize the ampleness for a divisor on a complex
  smooth projective variety with a universal model category
  associated with the divisor
  (\cite{lee-24-b}). 
\end{remark}

%\bibliographystyle{amsalpha}
%\bibliography{../../../dat}

\providecommand{\bysame}{\leavevmode\hbox to3em{\hrulefill}\thinspace}
\providecommand{\MR}{\relax\ifhmode\unskip\space\fi MR }
% \MRhref is called by the amsart/book/proc definition of \MR.
\providecommand{\MRhref}[2]{%
  \href{http://www.ams.org/mathscinet-getitem?mr=#1}{#2}
}
\providecommand{\href}[2]{#2}

\end{document}